\documentclass[paper]{elsarticle}
\pdfoutput=1
\usepackage{color}
\usepackage{multirow}
\usepackage{graphicx}
\usepackage{subfigure}
\usepackage{epsfig}
\usepackage{latexsym}
\usepackage{rotating}
\usepackage{amsmath, amssymb}
\usepackage{pdflscape}
\usepackage{bbm}
\usepackage{afterpage}
\usepackage{tikz}
\usetikzlibrary{shapes,arrows}

\newtheorem{remark}{Remark}[section]

\usepackage{lineno,hyperref}
\modulolinenumbers[5]

\journal{Communications in Nonlinear Science and Numerical Simulation}

\begin{document}
\begin{frontmatter}
	\title{\bf First Passage and First Exit Times  for diffusion processes related to a general growth curve}	
	\author[1]{G. Albano\corref{cor1}%
              }
               \ead{pialbano@unisa.it}
    \author[2,5]{A. Barrera} \ead{antonio.barrera@uma.es}
	\author[3]{V. Giorno} \ead{giorno@unisa.it}
	\author[4,5]{P. Rom\'an-Rom\'an} \ead{proman@ugr.es}
	\author[4,5]{F. Torres-Ruiz} \ead{fdeasis@ugr.es}	
	\address[1] {Dipartimento di Studi Politici e Sociali, Universit$\grave{a}$ degli Studi di Salerno, Via Giovanni Paolo II n.\ 132, I-84084 Fisciano (SA), Italy}
    \cortext[cor1]{Corresponding author}
    \address[2] {Department of Mathematical Analysis, Statistics and Operations Research and Applied Mathematics,  University of M$\acute{a}$laga, Spain}
	\address[3] {Dipartimento di Informatica, Universit$\grave{a}$ degli Studi di Salerno, Via Giovanni Paolo II n.\ 132, I-84084 Fisciano (SA), Italy}
	\address[4] {Department of Statistics and Operations Research, Faculty of Sciences, University of Granada, 18071 Granada, Spain}
    \address[5] {Institute of Mathematics of the University of Granada (IMAG), Calle Ventanilla, 11, 18001, Granada, Spain }
	\begin{abstract}
Recently a general growth curve including the well known growth equations, such as Malthus, logistic, Bertallanfy, Gompertz, has been studied. We now propose two stochastic formulations of this growth equation. They are obtained starting from a suitable parametrization of the deterministic model, by adding an additive and multiplicative noise respectively. For these processes we focus attention on the First Passage Time from a barrier and on the First Exit Time from a region delimited by two barriers. We consider thresholds, generally time dependent, for which there exist closed-forms of the probability densities of the first passage time and of the first exit time.
	\end{abstract}
	\begin{keyword}
First-passage times, first-exit times, ordinary differential equations, growth curve.
	\end{keyword}
\end{frontmatter}
\section{Introduction}
The models for the description of growth phenomena, originally associated with the evolution of animal populations, currently play an important role in several fields such as economics, biology, medicine, ecology (see, for example, \cite{Banks_94},  \cite{Jarne_07}). For this reason numerous efforts are oriented to the development of progressively sophisticated mathematical models for the description of a particular type of behavior. Systematically considered growth curves are of exponential type and, among these, the logistic and Gompertz laws. Indeed, the population size must be characterized by a limit due to the carrying capacity that represents, in general terms, the limitation of the natural resources. These curves are monotonic and have a sigmoidal shape; this is due to the presence of an inflection point in which the curve changes from concave to convex. Logistic and Gompertz curves have a similar growth pattern, the main difference between them being the location of the inflection point: the Gompertz curve reaches this point in the first part of the growth cycle whereas the logistic curve reaches it at later times.

S-shaped time evolutions are also observed in dynamic economic phenomena, such as the diffusion of technological innovations and product life cycles. These dynamics are usually incorporated into formal models by ordinary differential equations of the type $\dot x=f(x)$, where $f(\cdot)$ satisfies suitable properties (Turner et al. \cite{Turner_76}).

{\color{red}
These curves do not always adequately represent the realistic patterns; for instance, during the last two centuries, economic growth happens in the form of irregular successive expansions and contractions, whose expansion phases are longer than its contractions. Therefore, various generalizations have been proposed  as the  curves introduced by Von Bertalanffy, by Richards or  as the Hyperlogistic and Blumberg curves (Mialik et al. \cite{Economic2015}).
The von Bertalanffy curve (\cite{Bertalanfy_38}) is mainly used for modeling both length and weight for some animal species as well as in the study of the evolution of tumors treated with radiovirotherapy (\cite{Richards3}). The Richards curve, which has also been used to model the growth of animals (see, for instance, K\"{o}hn et al. \cite{Richards1} and
Nahashon et al. \cite{Richards2}), it has been applied in epidemiology (see \cite{Richards4}).
}

The complexity and variety of phenomena that can be described through growth models make it necessary, on the one hand, to provide mathematical models capable of describing peculiarities of certain dynamical systems and, on the other hand, to find equations capable of describing phenomena common to the various models. \par

In this direction, recently, the general growth curve initially proposed in \cite{Turner_76} and later in \cite{Tsoularis} has been studied in \cite{Albano_2022}. The formulation includes several parameters, whose choice leads to a variety of models such as the classical cases of Malthusian, Richards, Gompertz, Logistic and some their generalizations. The study focuses on the effects of the involved parameters through both analytical
results and computational evaluations.  \par
{\color{red}
We point out that the existence of discrepancies between the proposed models and the observed data suggests the use of stochastic models, among which those associated with stochastic diffusion processes stand out. Related to the logistic curve, Tuckwell and Koziol \cite{Tuc87} show a summary of some diffusion processes, some of which are linked to specific applications like Demography (Artzrouni and Reneke \cite{Art90}), or energy consumption (Giovanis and Skiadas \cite{Ski99}). Also, Schurz \cite{Sch07} considers a more general version for the stochastic differential equation associated to logistic growth. In the case of the Bertalanffy curve, we can cite the works of Quiming and Pitchford \cite{Qim07} and Rom\'an et al. \cite{RomanRomeroTorres_10}. However, it is perhaps the Gompertz curve that has been the subject of the most in-depth studies. Indeed, the fact that the Gompertz curve is an excellent model for the description of tumor growth has motivated the introduction of several diffusion processes associated with it (see, for example, Lo \cite{Lo07} and Ferrante et al. \cite{Fer05}). On the basis of these processes, a great number of modifications have subsequently been made with the goal of describing the evolution of tumor growth in the presence of therapeutic treatments (see, for instance, Albano et al. \cite{Alb20} and references therein). On other hand, apart from the classic growth models, new curves have recently been introduced that show more flexibility in their behaviour. In this sense, Tabatabai et al. \cite{Tab05} have constructed the so-called hyperbolastic curves and shown their usefulness in the study of the evolution of tumor processes and stem cell growth (Tabatabai et al. \cite{Tab11}). Regarding the stochastic versions of these latest models, Barrera et al. \cite{Bar2020} present a joint vision of all of them from the perspective of stochastic differential equations.

We must note that the introduction of these type of models allows us to deepen the study of the dynamic phenomena under consideration, beyond the results that can be obtained through the deterministic versions. Indeed, the inclusion of a probability structure makes it possible to approach the inference of the models from observed sample data of the phenomenon. In particular, it is possible to estimate parametric functions that represent important characteristics associated with the curves under study, as well as confidence intervals for the estimates and predictions (see, for example Guti\'errez et al. \cite{Gut05}). In the same way, the fact of dealing with stochastic diffusion processes allows us to consider the study of temporal variables that represent concrete problems. Among these problems, and associated with growth phenomena, it is possible to calculate the distribution of the random variable that indicates the time in which the growth of the phenomenon reaches a certain value for the first time or, in the case of sigmoidal phenomena, the instant in which the inflection point is reached, which determines a change in growth behavior. These problems can be addressed by considering the above situations as first-passage times, following the lines drawn by Buonocore et al.\cite{Buo87} in the homogeneous case and by Guti\'errez et al. \cite{Gut95} in the inhomogeneous one.
}

In the present paper,  following the approach in \cite{Barrera2020}, we proposed two different stochastic models based on the deterministic curve  studied in \cite{Albano_2022} obtained by introducing a multiplicative or an additive noises in the growth equation. In particular, the two stochastic models obtained are of the diffusive type and are characterized by the same mean which coincides with the deterministic curve. Also, they are lognormally and normally distributed, respectively; therefore the main probabilistic characteristics, such as transition probability density function and the related conditional moments, can be derived.\par

For the obtained stochastic processes, we analyze the first passage time (FPT) through a time dependent boundary and the first exit time (FET) from a  region bordered by two time dependent boundaries. This problem is particularly relevant in applications in which it is of interest to know the time for the population size to reach a fixed level (see, for example, \cite{MathematicsLoan}). Furthermore, in biological and economic contexts, it is often useful to identify the time in which the population first reaches a value coinciding with a percentage of its average size.\\

The paper is organized as follows. In Section 2 the deterministic model is described and a new parametrization is provided in order to prepare the ground for the stochastic extension. Since the obtained stochastic processes are related to Gauss Markov (GM) processes, in Section 3 a brief review of the GM processes is made in which the necessary preliminary background and notation is provided. In Section 4, by including an additive noise to the deterministic growth equation, a lognormal process is obtained and for it the FPT and the FET is anlysed providing suitable boundaries for wich the pdf's have closed-forms. In Section 5, by means of a multiplicative noise, an Ornstein-Uhlenbeck process is obtained and the FPT and FET analysed. A numerical anlysis is provided in Section 6, while some concluding remarks close the paper in Section 7.
\section{The deterministic model}
We consider the general growth curve described in the papers \cite{Turner_76} and \cite{Albano_2022}:
\begin{equation}\label{ODE}
\dfrac{dx}{dt}=\gamma k^{n(p-1)}x^{1+n(1-p)}\left[1-\left(\dfrac{x}{k}\right)^n\right]^{p}, \quad x(t_0)=x_0,
\end{equation}
where $x_0$ denotes the population size at the initial time $t_0$, $k=\lim\limits_{t \to \infty}x(t)>0$ represents the
carrying capacity,  $\gamma$, $n$ and $p$ are shape-parameters subject to being  positive with $0<p<1+1/n$.\par
The solution of Eq. \eqref{ODE} is
\begin{equation}\label{sol1}
x(t)=\dfrac{k}{\left\{1+\left[\gamma n(p-1)(t-t_0)+A_n^{1-p}\right]^{\frac{1}{1-p}}\right\}^{1/n}}.
\end{equation}
where $A_n=\left(\dfrac{k}{x_0}\right)^{n}-1$ depends on the shape parameter $n$ and on the ratio between the
carrying capacity $k$ and the initial population size $x_0$. \par
In \cite{Albano_2022} the curve \eqref{ODE} was extensively studied, showing that it is able to generalize the most famous growth equations. In addition, it presents some peculiar behaviors for suitable choices of the parameters $n$ and $p$, such as indefinite growth in a finite time or a plateau in a time interval followed by an indefinite increase  or by a decrease  to the initial value $x_0$. In particular, this last behavior is interesting in the context of epidemics.\par
Encouraged by these interesting findings, in the following we consider a re-parametrization of \eqref{sol1} which lends itself better to random generalization, also in the light of the particular cases that the equation generalizes.\par
Specifically,  defining
\begin{equation}
\alpha=e^{-\gamma\,n},\qquad\eta=[A_n^{1-p}+n\,\gamma\,(1-p)\,t_0]^{1/(p-1)},\label{alpha_eta}
\end{equation}
the \eqref{sol1} can be rewritten in the following form:
\begin{equation}
x(t)=x_0\,\frac{g(t_0)}{g(t)},\label{different_param}
\end{equation}
with
\begin{equation}
g(t)=\bigl\{\eta+[1+\eta^{1-p}\,\ln\alpha\,(1-p)\,t]^{1/(1-p)}\bigr\}^{1/n}.\label{gg}
\end{equation}
%
%
We note that Eq.~\eqref{sol1} comes from \eqref{different_param}. Indeed, from \eqref{alpha_eta} and \eqref{gg}, one has:
\begin{eqnarray}
&&g^n(t)=\Bigl[\frac{1}{A_n^{1-p}+n\,\gamma\,(1-p)\,t_0}\Bigr]^{1/(1-p)}+\Bigl[1-\frac{\gamma\,n\,(1-p)\,t}{A_n^{1-p}+n\
\gamma\,(1-p)\,t_0}\Bigr]^{1/(1-p)},\nonumber\\
&&\hspace{0.8 cm}=\frac{1+[A_n^{1-p}-\gamma\,n\,(1-p)\,(t-t_0)]^{1/(1-p)}}{[A_n^{1-p}+n\,\gamma\,(1-p)\,t_0]^{1/(1-p)}}
\end{eqnarray}
so that
\begin{equation}
\frac{g(t_0)}{g(t)}=\Bigl\{\frac{A_n+1}{1+[A_n^{1-p}-\gamma\,n\,(1-p)\,(t-t_0)]^{1/(1-p)}}\Bigr\}^{1/n}. \label{rem1_eq}
\end{equation}
Finally, since $A_n+1=\left(\dfrac{k}{x_0}\right)^{n}$,  from \eqref{different_param} and
\eqref{rem1_eq} Eq.~\eqref{sol1} follows.

We note that the proposed parametrization is coherent with that one provided in \cite{Roman_2015}  for the Richards curve (for $p\to 1$) and logistic curve (for $n\to 1$ and $p\to 1$).
\begin{remark}
Equation~\eqref{different_param} satisfies the following ordinary differential equation
\begin{equation}
\frac{dx}{dt}=h(t)\, x(t),\label{New_ODE}
\end{equation}
where
\begin{eqnarray}
&&h(t)=-\frac{g'(t)}{g(t)}=-\frac{\eta^{1-p}\,\ln\alpha\,[1+\eta^{1-p}\,\ln\alpha\,(1-p)\,t]^{p/(1-p)}}{n\,[g(t)]^n}\nonumber\\
&&\hspace{0.8cm}=-\frac{d}{dt}\ln g(t).\label{h}
\end{eqnarray}
\end{remark}
Indeed, by deriving \eqref{different_param} to respect to $t$, we obtain
$$
\frac{dx}{dt}=-x_0 g(t_0)\frac{g'(t)}{g^2(t)},
$$
and hence \eqref{New_ODE} holds.

Equation \eqref{New_ODE}  represents a Malthus growth with time dependent  fertility.
In the following two  stochastic extensions of the deterministic growth curve given in
\eqref{different_param} will be addressed. The idea is to consider two diffusion processes whose means are equal to \eqref{New_ODE}. This allows to reproduce different patterns starting from real data. \par
The starting point is \eqref{New_ODE} to which we add a multiplicative noise and an additive one, respectively. The resulting processes have different characteristics. As we will see the first process is a version of a lognormal diffusion process with exogenous factors, while the second one is a time inhomogeneous Ornstein-Uhlenbeck (OU) process. For both the processes we analyze the first passage time through a time dependent boundary.
We remark that the obtained stochastic processes are related to Gauss Markov (GM) processes, in the sense that the OU process is a GM process with a degenerate initial distribution, while the lognormal can be transformed into a GM process. This link permits to use all the techniques known in the literature for GM processes. With this in mind, a short review of these processes is addressed in the next section, in which the necessary preliminary background and notation is provided. \\
\section{Gauss-Markov processes in a nutshell}
Let $\{X(t), t\in T\}$, with $T$ denoting a continuous parameters set, be a continuous Gauss-Markov (GM) process such that the following properties hold:
\begin{itemize}
\item $m(t)=E[X(t)]$ is continuous for $t\in T$,
\item the covariance function $c(s,t)= E\{[X(s)-m(s)][X(t)-m(t)]\}$ is a continuous function for $(s,t)\in T^2$,
\item $\{X(t)\} $ is non-singular except possibly at the end points of $T$.
\end{itemize}
This last point means that if $T=[a,b]$ then $X(t)$ has a non-singular normal distribution except possibly for the points $a$ and $b$ where $X(t)$ could be degenerate in $m(t)$. \\
The Gauss-Markov processes satisfy well known properties (see, for example, \cite{Doob1949} and \cite{Mehr1965}) that will be used in the following. In particular,
\begin{itemize}
\item a Gaussian process is Markovian if and only if the covariance function can be expressed as
\begin{equation}
c(s,t)=k_1(s)\,k_2(t), \qquad s\leq t\label{covariance}
\end{equation}
where $k_1(t)$ and $k_2(t)$ are such that
\begin{equation}
r(t)=\frac{k_1(t)}{k_2(t)}\label{r}
\end{equation}
is a monotonically increasing function by virtue of the Cauchy-Schwarz inequality, with $k_1(t)\,k_2(t)>0$ because of the assumed non-singularity of the process in the interior of $T$;
\item the transition pdf $f(x,t|y,\tau)$ of a Gauss-Markov process is a normal density with conditional mean and variance:
\begin{eqnarray}
&&E[X(t)|X(\tau)=y]=m(t)+\frac{k_2(t)}{k_2(\tau)}\,[y-m(\tau)],\nonumber\\
&&Var[X(t)|X(\tau)=y]=k_2(t)\Bigl[k_1(t)-\frac{k_2(t)}{k_2(\tau)}\,k_1(\tau)\Bigr],
\end{eqnarray}
respectively.
\end{itemize}
%
 %
We note that the transition pdf $f(x,t|y,\tau)$  of a Gauss-Markov process satisfies the Fokker Plank equation (see \cite{Nobile_2001}):
$$
\frac{\partial f(x,t|y,\tau)}{\partial t}=-\frac{\partial}{\partial x} \bigl[B_1(x,t) f(x,t|y,\tau)\bigr]+\frac{1}{2}\frac{\partial^2}{\partial x^2} \bigl[B_2(t) f(x,t|y,\tau)\bigr]
$$
with the associated initial condition
$\lim_{\tau\to t}f(x,t|y,\tau)=\delta(x-y)$, where $B_1(x,t)$ and $B_2(t)$ are the drift and the infinitesimal variance of the process and they are given by:
\begin{equation}
B_1(x,t)= m'(t)+[x-m(t)]\frac {k'_2(t)}{k_2(t)},\qquad B_2(t)=k_2^2(t)\,r'(t).\label{mom_Gauss_Markov}
\end{equation}
In the following subsections we provide a  brief overview of the problems related to the FPT and to the FET  for the GM processes.
\subsection{FPT for GM processes}\label{sub3.1}
Let $s(t)$ be a continuous function and let $X(t_0)=x_0\neq s(t_0)$. The random variable FPT of $X(t)$ through the boundary $s(t)$
is {\color{red} defined as}
$$
T_{x_0,s(t)} =\left\{\begin{array}{ll}
\displaystyle {\inf_{t\geq t_0}\{t\;:\; X(t)>s(t)|X(t_0)=x_0\}},&x_0<s(t_0)\\
\displaystyle {\inf_{t\geq t_0}\{t\;:\; X(t)<s(t)|X(t_0)=x_0\}},&x_0>s(t_0)
\end{array}
\right.
$$
and $g[s(t),t|x_0,t_0]$ is the FPT pdf. \\
Following \cite{Nobile_2001}, we have the following remark.
\begin{remark}\label{remark3}
If $s(t), m(t), k_1(t)$ and $k_2(t)$ are $C^1(T)$-class, for $x_0<s(t_0)$, then $g[s(t),t|x_0,t_0]$ is solution of the following second kind Volterra integral equation:
\begin{equation}
\label{g}
g[s(t),t|x_0,t_0]=-2\Psi[s(t),t|x_0,t_0]+2\int_{t_0}^t g[s(\tau),\tau|x_0,t_0]\,\Psi[s(t),t|s(\tau),\tau]\,d\tau,
\end{equation}
where
\begin{eqnarray}
&&\Psi[s(t),t|y,\tau]=\Bigl\{\frac{s'(t)-m'(t)}{2}-\frac{s(t)-m(t)}{2}
\frac{k'_1(t)k_2(\tau)-k'_2(t)k_1(\tau)}{k_1(t)k_2(\tau)-k_2(t)k_1(\tau)}\nonumber\\
&&\hspace{0.8cm}-\frac{y-m(\tau)}{2}\frac{k'_2(t)k_1(t)-k_2(t)k'_1(t)}{k_1(t)k_2(\tau)-k_2(t)k_1(\tau)}\Bigr\} f[s(t),t|y,\tau].
\label{psi}
\end{eqnarray}
Further, $\Psi[s(t),t|s(\tau),\tau]=0$ for $\tau,t\in T$ with $\tau<t$ if and only the threshold $s(t)$ is chosen in the following way:
\begin{equation}
\widetilde s(t)=m(t)+d_1k_1(t)+d_2k_2(t),\qquad d_1,d_2\in\mathbb{R}.\label{s(t)}
\end{equation}
\end{remark}
Hence, from \eqref{g}, by choosing  $s(t)$ as in \eqref{s(t)}, one has the following closed-form for the FPT density:
\begin{equation}
g[\widetilde s(t),t|x_0,t_0]=\frac{\widetilde s(t_0)-x_0}{r(t)-r(t_0)}\frac{k_2(t)}{k_2(t_0)}r'(t)f[\widetilde s(t),t|x_0,t_0], \qquad x_0<\widetilde s(t_0),\label{closed_form_g}
\end{equation}
where $r'(t)$ is the derivative of the function $r(t)$ defined in \eqref{r} and $f$ is the transition pdf of $X(t)$.
%
\subsection{FET problem for GM processes}\label{sub3.2}
Let $s_1(t)$ and  $s_2(t)$ be continuous functions  such that
$s_1(t)<s_2(t)$ for $t>t_0$. Assuming that $X(t_0)=x_0$ is such that $s_1(t_0)<x_0<s_2(t_0)$, we analyze the evolution of
the process $X(t)$ defined in \eqref{mom_Gauss_Markov} in the presence of two absorbing boundaries in $s_1(t)$ and  $s_2(t)$. Specifically, we consider the following FET random variables:
\begin{eqnarray}
&&\hspace{-0.7cm}T_1=\inf_{t\geq t_0}\{t:X(t)< s_1(t);\
X(\theta)< s_2(\theta),\forall\ \theta \in (t_0, t)\}
\quad{\rm (FET\;from\; above)}\nonumber\\
&&\nonumber\\
&&\hspace{-0.7cm}T_2=\inf_{t\geq t_0}\{t : X(t)> s_2(t);\
X(\theta)>s_1(\theta),\forall\ \theta \in (t_0, t)\}
\quad{(\rm FET\;from\; below)}\nonumber\\
&&\nonumber\\
&&\hspace{-0.7cm}T_{1,2}=\inf \{{T_1,\ T_2}\}\hspace{0.9cm}{(\rm FET)},\nonumber
\end{eqnarray}
characterized by pdf's
\begin{eqnarray}
\label{def_gamma_L_AB}&& \gamma_1(t|x_0,t_0)=\frac{\partial}{\partial t}\
P(T_1< t),\qquad\gamma_2(t| x_0,t_0)=\frac{\partial}{\partial t}\
P(T_2< t)\nonumber \\
&&\gamma(t|x_0,t_0)=\gamma_1(t|x_0,t_0)+\gamma_2(t|x_0,t_0),
\end{eqnarray}
respectively.
{\color{red} Clearly, the functions  $\gamma_i$  depend on the barriers $s_1(t)$ and $s_2(t)$  although, for simplicity, we have omitted this dependence in the notation.} \par\noindent
Following \cite{Nobile_2006}, we have the following remark.
%
\begin{remark}\label{remark4}
If $s_1(t), s_2(t), m(t), k_1(t)$ and $k_2(t)$ are $C^1(T)$-class functions and $s_1(t)<s_2(t)$ for all $t \in T$, with $s_1(t_0)<x_0<s_2(t_0)$, the functions $\gamma_1(t|x_0,t_0)$ and $\gamma_2(t|x_0,t_0)$ are solutions of the following second kind Volterra integral equations:
\begin{eqnarray}
\label{gamma1}
&&\hspace{-1.2cm}\gamma_1(t|x_0,t_0)=2\Psi_1(t|x_0,t_0)\nonumber\\
&&-2\int_{t_0}^t \bigl\{\gamma_1(\tau|x_0,t_0]\,\Psi_1[t|s_1(\tau),\tau]+
\gamma_2(\tau|x_0,t_0]\,\Psi_1[t|s_2(\tau),\tau]\bigr\}\,d\tau,\nonumber\\
&&\hspace{-1.2cm}\gamma_2(t|x_0,t_0)=-2\Psi_2(t|x_0,t_0)\nonumber\\
&&+2\int_{t_0}^t \bigl\{\gamma_1(\tau|x_0,t_0]\,\Psi_2[t|s_1(\tau),\tau]+
\gamma_2(\tau|x_0,t_0]\,\Psi_2[t|s_2(\tau),\tau]\bigr\}\,d\tau,\nonumber
\end{eqnarray}
where, for $ j=1,2$, one has:
\begin{eqnarray}
&&\hspace{-1.2cm}\Psi_j[t|y,\tau]=\Bigl\{\frac{s'_j(t)-m'(t)}{2}-\frac{s_j(t)-m(t)}{2}
\frac{k'_1(t)k_2(\tau)-k'_2(t)k_1(\tau)}{k_1(t)k_2(\tau)-k_2(t)k_1(\tau)}\nonumber\\
&&-\frac{y-m(\tau)}{2}\frac{k'_2(t)k_1(t)-k_2(t)k'_1(t)}{k_1(t)k_2(\tau)-k_2(t)k_1(\tau)}\Bigr\} f[s_j(t),t|y,\tau].
\label{psij}
\end{eqnarray}
Further, if $x_0 = m(t_0) + a k_1(t_0) + c k_2(t_0)$ and if
$$
s_1(t) = m(t) + a k_1(t) + c_1 k_2(t), \quad s_2(t) = m(t) + a k_1(t) + c_2 k_2(t),
$$
with $s_1(t)<s_2(t)$, for all $t\in T$
then
\begin{eqnarray}\label{gamma}
&&\hspace{-1.2cm}\gamma(t|x_0, t_0) =\frac{k_2(t)}{r(t) - r(t_0)}r'(t)\sum_{n=-\infty}^\infty\exp\Bigl\{
-\frac{2 n^2 (c_2 - c_1)^2}{r(t) - r(t_0)}\Bigr\}\nonumber\\
&&\hspace{-0.6cm}\times \Bigl\{\Bigl[c - c_1 + 2n (c_2 - c_1)\Bigr]\exp\Bigl\{-\frac{2 n (c_2 - c_1) (c - c_1)}{r(t) - r(t_0)}\Bigr\}f[s_1(t), t | x_0, t_0]\nonumber\\
&&\hspace{-0.6cm}+\Bigl[c_2 - c- 2n (c_2 - c_1)\Bigr]\exp\Bigl\{\frac{2 n (c_2 - c_1) (c_2 - c)}{r(t) - r(t_0)}\Bigr\}f[s_2(t), t | x_0, t_0],\qquad
\end{eqnarray}
with $a, c, c_1, c_2 \in R$.
\end{remark}
The results given in Section \ref{sub3.1} and \ref{sub3.2} will be used in the following to analyze the FPT and the FET problems for the stochastic processes obtained adding some types of noise to the deterministic differential equation \eqref{different_param}.
\section{Multiplicative noise: Lognormal diffusion process}
Let $X_L(t)$ be the stochastic process obtained from the deterministic growth equation \eqref{different_param} by including a multiplicative noise. In particular, the introduction of  a  noise
with variance  $\sigma^2>0$  in the intrinsic fertility $h(t)$ leads to the following stochastic differential equation (SDE) (see, for instance, \cite{Albano2006}):
\begin{equation}
\label{SDEL}
dX_{L}(t)=h(t) X_{L}(t) dt+\sigma  X_{L}(t) dW(t), \qquad X_L(t_0)=x_0,
\end{equation}
where $W(t)$ is a standard Wiener process in $\mathbb R$, independent on {\color{red} the initial size of the population} and   $h(t)$ is the function defined in
\eqref{h}. {\color{red} We point out that the initial size is usually a random variable, however in many real applications the value $X_L(t_0)$ is known, so that we assume that $x_0$ is a degenerate random variable, i.e. known without errors.} From \eqref{SDEL} we deduce that $X_L(t)$ is a diffusion process with state space $[0,+\infty)$, having drift and infinitesimal variance:
\begin{equation}
A_1(x,t)=h(t)\,x,\qquad A_2(x)=\sigma^2\,x^2.\label{log_process}
\end{equation}
The solution of \eqref{SDEL} is a inhomogeneous diffusion  lognormal process:
\begin{eqnarray}
&&X_{L}(t)=x_0\, \exp\Biggl\{ \int_{t_0}^t h(\xi) d\xi-\frac{\sigma^2}{2}(t-t_0)+\sigma\bigl[W(t)-W(t_0)\bigr]\Biggr\}
\nonumber\\
&&\hspace{1cm}=x_0\, \frac{g(t_0)}{g(t)}\exp\Biggl\{-\frac{\sigma^2}{2}(t-t_0)+\sigma\bigl[W(t)-W(t_0)\bigr]\Biggr\},
\end{eqnarray}
where \eqref{h} has been used.\\
The process $X_L(t)$ can be reduced to a Wiener process $Z(t)$ characterized by drift and infinitesimal variance
\begin{equation}
 B_1=0,\qquad B_2=\sigma^2,
\label{Wiener_process}
 \end{equation}
by means of Ito's lemma  and the following transformation
\begin{eqnarray}
&& z=\log x-\int^th(\xi)\,d\xi+\frac{\sigma^2}{2}\,t.
\label{transformation}
\end{eqnarray}
 The process $Z(t)$ is a GM process with $B_1(x,t)=0$ and $B_2(t)=\sigma^2 $, so that from Remark \eqref{mom_Gauss_Markov}, we have:
$$
k_1(t)=\sigma^2 t,\qquad k_2(t)=1,\qquad m(t)=0.
$$
The transition pdf of $X_L(t)$, denoted by $f_L(x,t|y,\tau)=\frac{d}{dx}P[X_L(t)<x\mid X_L(\tau)=y]$, is a lognormal density. Specifically, for all $t>\tau\geq t_0$
the conditional random variable $X_L(t)\mid X_L(\tau)=y$ follows the lognormal distribution  $\Lambda_1\bigl(M_L(t|\ln y, \tau), \sigma^2(t-\tau)\bigr)$:
\begin{equation}
f_L(x,t|y,\tau)=\frac{1}{x\sqrt{2 \pi \sigma^2(t-\tau)}}\exp\Biggl\{-\frac{[\ln x-M_L(t|\ln y,\tau)]^2}{2 \sigma^2(t-\tau)}\Biggr\},
\qquad x,y\in\mathbb{R}^+,\label{f_L}
\end{equation}
with
$$
M_L(t|\ln y,\tau)= \ln y- \frac{\sigma^2}{2}(t-\tau)+\int_\tau^t h(\xi) d\xi=  \ln y+\ln \frac{g(\tau)}{g(t)}- \frac{\sigma^2}{2}(t-\tau).
$$
Moreover, the conditional cumulative transition distribution of $X_L(t)$ is given by
$$
F_L(x,t|y,\tau)=\int_0^x f_L(z,t|y,\tau)\,dz =\frac{1}{2}\biggl\{1+{\rm Erf}
\biggl[\frac{\ln x-M_L(t|\ln y,\tau)}{\sqrt{2\,\sigma^2(t-\tau)}}\biggr]\biggr\},
$$
where ${\rm Erf}(z)= \frac{2}{\sqrt\pi}\int_0^ze^{-\xi^2}\, d\xi$ is the error function.
The conditional moments for $X_L(t)$ are:
$$
\mathbb{E}[X^n_L(t)|X(\tau)=y]=\exp\bigl\{n\,M_L(t|\ln y,\tau)+\frac{n^2}{2}\,\sigma^2(t-\tau)\bigr\},
$$
from which we can easily obtain the conditional mean and  variance:
\begin{eqnarray}
&&\hspace{-0.8cm}\mathbb{E}[X_L(t)|X_L(\tau)=y]=y\,\frac{g(\tau)}{g(t)},\nonumber\\
&&\hspace{-0.8cm}{\rm Var}[X_L(t)|X_L(\tau)=y]=\Bigl[y\,\frac{g(\tau)}{g(t)}\Bigr] ^2\Bigl[\exp\{\sigma^2(t-\tau)\}-1\Bigr].\nonumber
\end{eqnarray}
We point out that the conditional mean has the same trend of the solution \eqref{different_param} of the deterministic equation.
 \subsection{First passage time problem}
Let $T_{x_0,s(t)}^L $ be the random variable FPT of $X_L(t)$ through a continuous time-varying boundary $s(t)$
and let $g_L[s(t),t|x_0,t_0]$ be the FPT pdf. We can study the FPT problem for the process $X_L(t)$ starting from the results concerning the FPT of the transformed Wiener process $Z(t)$.
In particular,  from \eqref{s(t)} and \eqref{closed_form_g}, by choosing  $\widetilde s(t)=a\,t+b$, one has:
$$
g_Z[\widetilde  s(t),t|\widetilde y,\tau]=\frac{|\widetilde s(\tau)-\widetilde y|}{t-\tau}f_Z[\widetilde s(t),t|\widetilde y,\tau], \qquad \widetilde y\neq \widetilde s(\tau)
$$
being $f_Z[\widetilde s(t),t|\widetilde y,\tau]$ the transition pdf of the Wiener process $Z(t)$.
Then, recalling the transformation \eqref{transformation}, we can obtain the FPT density of $X_L(t)$ through the boundary
\begin{equation}
s(t)=A\exp\Bigl\{Bt+\int_0^t h(\xi)\,d\xi\Bigr\},\qquad A>0,\;B\in {\mathbb R},\label{soglia_L}
\end{equation}
with $h(t)$ given in \eqref{h}. In particular we have
\begin{eqnarray}
&&\hspace{-2.0cm}g_L[s(t),t|x_0,t_0]={\Big|\ln {s(t_0)\over x_0}\Bigr|\over\sqrt{2\pi\sigma^2(t\!-\!t_0)^3}}\nonumber\\
&&\hspace{-.0cm}\times \exp\biggl\{-{\Bigl[({\sigma^2\over 2}\!+\!B)(t\!-\!t_0)\!+\!\ln {s(t_0)\over x_0}\Bigr]^2\over2\sigma^2(t\!-\!t_0)}\biggr\}, \;\; s(t_0)\neq x_0.\label{closed_g_L}
\end{eqnarray}
Moreover,  by choosing  in (\ref{soglia_L})
$B=0$ and $ A=\nu\,x_0\exp\Bigl\{-\int_0^{t_0}h(\xi)\,d\xi\}\Bigr\}=\nu\,x_0\frac{g(0)}{g(t_0)}$, with $g(t)$ defined in \eqref{gg},
one has
\begin{equation}
s(t) = \nu E[X_L(t)|X_L(t_0)=x_0]=\nu\, x_0\,\exp\Bigl\{\int_{t_0}^t h(\xi)\,d\xi\Bigr\}=\nu\, x_0 \frac{g(t_0)}{g(t)},
\label{boundary}
\end{equation}
that, for $0<\nu<1$, represents a percentage of the conditional mean of the process $X_L(t)$. Therefore, for the process $\{X_L(t);\,t\geq t_0\}$ characterized  by infinitesimal moments \eqref{log_process},
the FPT pdf through  the boundary (\ref{boundary})  is given by
\begin{equation}
g_L[s(t),t|x_0,t_0]={|\ln \nu |\over\sqrt{2\pi\sigma^2(t-t_0)^3}}
\exp\Bigl\{-{[\sigma^2\,(t-t_0)/2+\ln \nu]^2\over2\sigma^2(t-t_0)}\Bigr\}.\label{proportional}
\end{equation}
We note that \eqref{closed_g_L} and \eqref{proportional} can be obtained also following alternative procedure as those one proposed in  \cite{Sato_1989}, \cite{Gut95}. In these papers, such as in \cite{Nobile_2001},  one can also find procedures  to obtain good numerical approximations to the FPT density.
%
%
\subsection{First exit time problem}
Let $T_1^L, T_2^L$ and $T_{1,2}^L$ be the FET random variables of the process $X_L(t)$ through the continuous boundaries $s_1(t)$ and  $s_2(t)$ $(s_1(t)<s_2(t)$ for $t>t_0)$.
Since the process $X_L(t)$ is
time-inhomogeneous, procedures to analyze the FET problem for it do not appear to be present in literature. However, the problem can be solved  recalling that  the  logarithmic  transformation \eqref{transformation} applied to $X_L(t)$ leads to a Wiener process $Z(t)$ with infinitesimal moments \eqref{Wiener_process}. \\
Following the results provided in Section 3.2, we can obtain closed-forms for the FET pdf $\gamma_Z$ by choosing suitable boundaries as described in the following.\\
From Remark \ref{remark4}, setting
\begin{equation}
\widetilde s_i(t)=c_i+\alpha t\qquad i=1,2;\label{boundaries_Z}
\end{equation}
and $z_0=\alpha t_0+c$
with $c_1<c<c_2,$ the FET pdf $\gamma^Z(t|z_0,t_0)=\gamma_1^Z(t|z_0,t_0)+\gamma_2^Z(t|z_0,t_0)$ of $Z(t)$
from the interval $\bigl(\widetilde s_1 (t), \widetilde s_2(t)\bigr)$ is:
\begingroup
\footnotesize
\begin{eqnarray}
&&\hspace{-0.8cm}\gamma^Z(t|z_0,t_0)={1\over \sqrt{2\,\pi\sigma^2(t-t_0)^3}}\sum_{n=-\infty}^{+\infty}
\exp\Bigl\{-{2n^2(c_2-c_1)^2\over \sigma^2(t-t_0)}\Bigr\}\nonumber\\
&&\times\biggl\{[c-c_1+2n(c_2-c_1)]\,\exp\Bigl\{-{2n(c_2-c_1)(c-c_1)\over \sigma^2(t-t_0)}
\Bigr\}\exp\Bigl\{-{[\alpha(t-t_0)+c_1-c]^2\over 2\sigma^2(t-t_0)}\Bigr\}\nonumber\\
&&\hspace{0.4cm}+
[c_2-c-2n(c_2-c_1)]\,
\exp\Bigl\{{2n(c_2-c_1)(c_2-c)\over \sigma^2(t-t_0)}\Bigr\} \exp\Bigl\{-{[\alpha(t-t_0)+c_2-c]^2\over 2\sigma^2(t-t_0)}\Bigr\}
\biggr\}.\nonumber\\
\label{gamma_lineare}
\end{eqnarray}
\endgroup
For some choices of the parameters Eq.~\eqref{gamma_lineare} can be simplified; this happens, for example, if the boundaries $\widetilde s_i(t)$  are time independent  and $z_0$ is the midpoint of the interval $(\widetilde s_1,\widetilde s_2)$. In other words, if $\alpha=0$ and $z_0=(c_1+c_2)/2$ then
$\gamma_1^Z(t,|\widetilde y,\tau)=\gamma_2^Z(t,|\widetilde y,\tau)=\gamma^Z(t|\widetilde y,\tau)/2$ so from \eqref{gamma_lineare} we obtain:
\begingroup
\footnotesize
\begin{eqnarray}
&&\hspace{-0.8cm}\gamma^Z(t|z_0,t_0)=
{c_2-c_1\over \sqrt{2\,\pi\sigma^2(t-t_0)^3}}\exp\Bigl\{-{[c_2-c_1]^2\over 8\sigma^2(t-t_0)}\Bigr\}\Biggl\{1+\sum_{n=1}^{+\infty}
\exp\Bigl\{-{2n^2(c_2-c_1)^2\over \sigma^2(t-t_0)}\Bigr\}\nonumber\\
&&\times\biggl[[1+4n(c_2-c_1)]\,\exp\Bigl\{-{n(c_2-c_1)^2\over \sigma^2(t-t_0)}
\Bigr\}
+
[1-4n(c_2-c_1)]\,
\exp\Bigl\{{n(c_2-c_1)^2\over \sigma^2(t-t_0)}\Bigr\} \biggr]\Biggr\}.\nonumber
\label{gamma_lineare_centrale}
\end{eqnarray}
\endgroup
Making use of the results just shown, we look at the FET problem of the process $X_L(t)$ through the boundaries $s_i(t)=\exp\Bigl\{c_i+\int^th(\xi)\,d\xi+\alpha t-\frac{\sigma^2}{2}t\Bigr\}$ for $i=1,2$ obtained from \eqref{transformation}  and \eqref{boundaries_Z}. We note that the boundaries $s_i(t)$ have the same functional form of the  \eqref{soglia_L} since they can be expressed as
\begin{equation}
s_i(t)=A_i\,\exp\Bigl\{B\,t+\int_0^t h(\xi)\,d\xi\Bigr\}\label{boundaries_L}
\end{equation}
 with $A_i$ related to $c_i$, whereas the constant $B$ is expressible in terms of  the constants $\alpha$ and of the intensity of the noise $\sigma^2$. As in the case of a single boundary, we can choice $A_i$ and $B$ such that $S_i(t)$ become percentages of the conditional mean of the process. Indeed, if $B\equiv \alpha-\sigma^2/2=0$ and $ A_i=\nu_i\,x_0\exp\Bigl\{-\int_0^{t_0}h(\xi)\,d\xi\}\Bigr\}$, for $i=1,2$ one has
\begin{equation}
s_i(t) = \nu _i E[X_L(t)|X_L(t_0)=x_0]=\nu_i\, x_0\,\int_{t_0}^t h(\xi)\,d\xi=\nu_i\, x_0\,\frac {g(t_0)}{g(t)}\qquad \nu_1<\nu_2.\label{boundaries_L_mean}
\end{equation}
Therefore, from \eqref{gamma_lineare} by choosing $\alpha=\sigma^2/2$ and $c_i=\ln\nu_i+\ln x_0-\int^{t_0} h(\xi)\,d\xi$ and $x_0=\exp\Bigl\{\nu\int_0^{t_0} h(\xi)\,d\xi\Bigr\}=\frac{g(0)}{g(t_0)} \rm{e}^{\nu}$, we obtain the FET pdf of the process $X_L(t)$ characterized  by infinitesimal moments \eqref{log_process} through  boundaries \eqref{boundaries_L_mean}:
\begingroup
\footnotesize
\begin{eqnarray}
&&\hspace{-1.3cm}\gamma_L(t|x_0,t_0)={1\over \sqrt{2\,\pi\sigma^2(t-t_0)^3}}\sum_{n=-\infty}^{+\infty}
\exp\Bigl\{-{2n^2(\ln\frac{\nu_2}{\nu_1})^2\over \sigma^2(t-t_0)}\Bigr\}\nonumber\\
&&\hspace{-0.3cm}\times\biggl\{\Bigl[\ln\frac{\nu}{\nu_1}+2n\ln\frac{\nu_2}{\nu_1}\Bigr]\,\exp\Bigl\{-\frac{2n\ln \frac{\nu_2}{\nu_1}\ln \frac{\nu}{\nu_1}}{\sigma^2(t-t_0)}\Bigr\}\exp\Bigl\{-\frac{\bigl[\frac{\sigma^2}{2}(t-t_0)+\ln\frac{\nu_1}{\nu}\bigr]^2}{2\sigma^2(t-t_0)}\Bigr\}\nonumber\\
&&\hspace{-0.3cm}+\biggl\{\Bigl[\ln\frac{\nu_2}{\nu}-2n\ln\frac{\nu_2}{\nu_1}\Bigr]\,\exp\Bigl\{\frac{2n\ln \frac{\nu_2}{\nu_1}\ln \frac{\nu_2}{\nu}}{\sigma^2(t-t_0)}\Bigr\}\exp\Bigl\{-\frac{\bigl[\frac{\sigma^2}{2}(t-t_0)+\ln\frac{\nu_2}{\nu}\bigr]^2}{2\sigma^2(t-t_0)}\Bigr\}\biggr\}.\quad
\label{gamma_percentuali}
\end{eqnarray}
\endgroup
%
\section{Additive noise: Ornstein-Uhlenbeck diffusion process}
Let $X_G(t)$ be the stochastic process obtained from the deterministic growth equation \eqref{different_param} by including an additive noise. Specifically,
starting from Eq.~\eqref{New_ODE}, we obtain a stochastic generalization by introducing a white noise with variance $\sigma^2$, where $\sigma>0$ represents the width of environment fluctuations.  In this way, we have the following SDE:
\begin{equation}
\label{SDEG}
dX_{G}(t)=h(t) X_{G}(t) dt+\sigma  dW(t), \qquad X_G(t_0)=x_0,
\end{equation}
where $W(t)$ is a standard Wiener process in $\mathbb R$, independent by $x$ and the function  $h(t)$ is defined in
\eqref{h}. From \eqref{SDEG} we conclude that $X_G(t)$ is a diffusion process with state space $\mathbb{R}$, and infinitesimal moments:
\begin{equation}
B_1(x,t)=h(t)\,x,\qquad B_2(x)=\sigma^2.\label{gauss_process}
\end{equation}
Making use of the transformation
$$
Z(t)= \exp\Bigl\{-\int^t h(\vartheta)\, d\vartheta\Bigr\}\,X_G(t)
$$
from \eqref{SDEG} we obtain
\begin{equation}
\label{SDEG_trasformata}
dZ(t)=\sigma  \exp\Bigl\{-\int^t h(\vartheta) \,d\vartheta\Bigr\} dW(t),
\end{equation}
so that
\begin{equation}
Z(t)=Z(t_0)+\sigma\int_{t_0}^t \exp\Bigl\{-\int^u h(\vartheta)\,d\vartheta\Bigr\}\, dW(u).
\end{equation}
Therefore, for the process $X_G(t)= \exp\Bigl\{\int^t h(\vartheta)\, d\vartheta\Bigr\}\,Z(t)$ we obtain:
\begin{equation}
X_G(t)=x_0\,\exp\Bigl\{\int^t_{t_0} h(\vartheta)\,d\vartheta\Bigr\}+\sigma\int_{t_0}^t \exp\Bigl\{\int_u^t h(\vartheta)\,d\vartheta\Bigr\}\, dW(u).
\end{equation}
The transition pdf of $X_G(t)$ is Gaussian, specifically, for $t_0\leq \tau<t$ the conditional variable
$X_G(t)\mid X_G(\tau)=y$  is distributed as   ${\cal N}\bigl(M_G(t| y,\tau), V_G(t|\tau)\bigr)$ and its pdf is
\begin{equation}
f_G(x,t|y,\tau)=\frac{1}{\sqrt{2 \pi V(t|\tau)}}\exp\Biggl\{-\frac{[x-M_G(t|y,\tau)]^2}{2 V(t|\tau)}\Biggr\},\qquad x,y\in\mathbb{R}\label{f_G}
\end{equation}
where
\begin{eqnarray}
&&\hspace{-0.9cm}M_G(t|y,\tau)=
y\,\exp\Bigl\{\int^t_{\tau} h(\vartheta)\,d\vartheta\Bigr\}=y\,\frac{g(\tau)}{g(t)}\nonumber\\
&&\hspace{-0.9cm}V_G(t|\tau)=
\sigma^2\int_\tau^t \exp\Bigl\{2\int_\tau^\vartheta h(u)\,du\Bigr\}\, d\vartheta=\sigma^2\int_\tau^t\Bigl[\frac{g(\tau)}{g(\vartheta)}\Bigr]^2\,d\vartheta.\label{moments_XG}
\end{eqnarray}
represent the conditional mean and the conditional variance of $X_G(t)$, respectively.
The conditional cumulative transition distribution is given by
$$
F_G(x,t|y,\tau)=\int_{-\infty}^x f_L(z,t|y,\tau)\,dz =\frac{1}{2}\biggl\{1+{\rm Erf}
\biggl[\frac{x-M_G(t| y,\tau)}{\sqrt{2\,V_G(t|\tau)}}\biggr]\biggr\},
$$
where ${\rm Erf}(z)$ is the error function. \par
In the following, we study the FPT and FET problems for $X_G(t)$. To this aim, we note that $X_G(t)$ is a Gauss-Markov process, so that we can study the FPT problem and the FET problem following the approach proposed in \cite{Nobile_2001} and \cite{Nobile_2006}, respectively. \par
\noindent Alternatively,  we can transform $X_G(t)$ into the GM process $Z(t)$ defined in \eqref{Wiener_process}.
\subsection{First passage time problem}
Let $T_{x_0,s(t)}^G $ be the random variable FPT of $X_G(t)$ through a boundary $s(t)$
and let $g_G[s(t),t|x_0,t_0]$ be the FPT pdf. \par
From
\eqref{gauss_process}, recalling \eqref{mom_Gauss_Markov}, we have:
$$
k_1(t)=\sigma^2k_2(t)\int^t\exp\Bigl\{-2\int^uh(\vartheta)\,d\vartheta\Bigr\}du=\sigma^2k_2(t)\int^t g^2(u)\,du
$$
$$
\hspace{-1.5cm}k_2(t)=\exp\Bigl\{\int^t h(\vartheta)\,d\vartheta\Bigr\}=\frac{1}{g(t)}, \qquad m(t)=0.
$$
A closed-form of the FPT pdf of $X_G(t)$ can be obtained making use of Remark \ref{remark3}. In particular, by choosing the boundary
\begin{eqnarray}
&&\hspace{-1.5cm}s(t)=\exp\Bigl\{\int^t h(\vartheta)\,d\vartheta\Bigr\}\biggl[A+B\sigma^2\int^t\exp\Bigl\{-2\int^uh(\vartheta)\,d\vartheta\Bigr\}du\biggr]\nonumber\\
&&=\frac{1}{g(t)}\Bigl\{A+ B\sigma^2 \int^t g^2(u)du\Bigr\},
\label{boundary_XG}
\end{eqnarray}
and by using \eqref{r},
the FPT pdf of $X_G(t)$ through $s(t)$ is given by
\begin{eqnarray}
&&\hspace{-1.5cm}g_G[s(t),t|x_0,t_0]=\frac{|s(t_0)-x_0|}{r(t)-r(t_0)}\,\frac{k_2(t)}{k_2(t_0)}\, r'(t)\,f_G[s(t),t|x_0,t_0]\nonumber\\
&&\hspace{1.2cm}=[g(t)]^2\frac{|s(t_0)-x_0|}{\int_{t_0}^t [g(u)]^2\,du} f_G[s(t),t|x_0,t_0]
\label{gG}
\end{eqnarray}
with $f_G[s(t),t|x_0,t_0]$ given in \eqref{f_G}. \par
We note that the boundaries include a percentage of a conditional  mean of $X_G(t)$; indeed, by comparing  \eqref{moments_XG} and \eqref{boundary_XG} and by choosing $B=0$ and
$
A=\nu x_0 g(t_0)
$
we have
$s(t)=\nu M_G(t\mid x_0,t_0)$.
%
\subsection{First exit time problem}
Let $T_1^G, T_2^G$ and $T_{1,2}^G$ be the FET random variables of the process $X_G(t)$ through the continuous boundaries $s_1(t)$ and  $s_2(t)$ $(s_1(t)<s_2(t)$ for $t>t_0)$. \\
Following the results provided in Section 3.2, we can obtain closed forms for the FET pdf by choosing suitable boundaries.\\
Specifically, making use of Remark \ref{remark4} and by choosing
\begin{eqnarray}
&&\hspace{-0.5cm}s_i(t)=\exp\Bigl\{\int^t h(\vartheta)\,d\vartheta\Bigr\}\biggl[c_i+B\sigma^2\int^t\exp\Bigl\{-2\int^uh(\vartheta)\,d\vartheta\Bigr\}du\biggr]\nonumber\\
&&\hspace{0.3cm}=\frac{1}{g(t)}\Bigl\{c_i+ B\sigma^2 \int^t g^2(u)du\Bigr\}\qquad i=1,2
\label{Si_XG}
\end{eqnarray}
\begin{eqnarray}
&&x_0=\exp\Bigl\{\int^t h(\vartheta)\,d\vartheta\Bigr\}\biggl[c+B\sigma^2\int^t\exp\Bigl\{-2\int^uh(\vartheta)\,d\vartheta\Bigr\}du\biggr]\nonumber\\
&&\hspace{0.5cm}=\frac{1}{g(t)}\Bigl\{c+ B\sigma^2 \int^t g^2(u)du\Bigr\}\label{x0_G}
\end{eqnarray}
with $c_1<c<c_2,$ the FET pdf $\gamma_G(t|x_0,t_0)=\gamma_1^G(t|x_0,t_0)+\gamma_2^G(t|x_0,t_0)$ of $X_G(t)$
from the interval $\bigl(s_1 (t),  s_2(t)\bigr)$ is obtained from  \eqref{gamma}. In particular, $\gamma_G(t|z_0,t_0)$ is obtained from \eqref{gamma} by replacing  the transition pdf $f$ with $f_G$ given in \eqref{f_G} and recalling that for $X_G(t)$ one has $k_2(t)=\frac{1}{g(t)}$ ,
$$
r(t)=\sigma^2\int^t g^2(u)\,du,\qquad r^\prime(t)=\sigma^2 g^2(t).
$$
We point out that also in this case, we can identify the boundaries $s_i(t)$ as a percentages of the conditional average of the process given in \eqref{moments_XG}. Precisely, by choosing $B=0$ in \eqref{Si_XG}, we obtain
$$
s_i(t)=c_i x_0 \frac{g(t_0)}{g(t)}\qquad i=1,2.
$$
and, from \eqref{gamma}, recalling \eqref{f_G}, the FET pdf becomes:
\begingroup
\footnotesize
\begin{eqnarray}\label{gammaG}
&&\hspace{-1.1cm}\gamma_G(t|x_0, t_0) =
\frac{g(t)}{g(t_0)}\frac{1}{\sqrt{2\pi \sigma^2 \int_{t_0}^t \frac{1}{g^2(u)} du}}
\frac{1}{\int_{t_0}^t g^2(u) du}\sum_{n=-\infty}^\infty\exp\Bigl\{
-\frac{2 n^2 (c_2 - c_1)^2}{\sigma^2\int_{t_0}^t g^2(u) du}\Bigr\}\nonumber\\
&&\hspace{-0.6cm}\times \Biggl\{\Bigl[c - c_1 + 2n (c_2 - c_1)\Bigr]\exp\Bigl\{-\frac{2 n (c_2 - c_1) (c - c_1)}{\sigma^2\int_{t_0}^t g^2(u) du}\Bigr\}
\exp\Bigl\{-\frac{(c_1-1)^2 x_0^2}{2\sigma^2 g^2(t)\int_{t_0}^t \frac{1}{g^2(u)} du}\Bigr\}\nonumber\\
&&\hspace{-0.6cm}+\Bigl[c_2 - c- 2n (c_2 - c_1)\Bigr]\exp\Bigl\{\frac{2 n (c_2 - c_1) (c_2 - c)}{\sigma^2\int_{t_0}^t g^2(u) du}\Bigr\}\exp\Bigl\{-\frac{(c_2-1)^2 x_0^2}{2\sigma^2 g^2(t)\int_{t_0}^t \frac{1}{g^2(u)} du}\Bigr\}\Biggr\}.\nonumber\\
&&
\end{eqnarray}
\endgroup
%
\section{Numerical results}
\begin{figure}[htbp]
\centering
\subfigure[$p=1.5$]
{\includegraphics[width=0.45\textwidth]{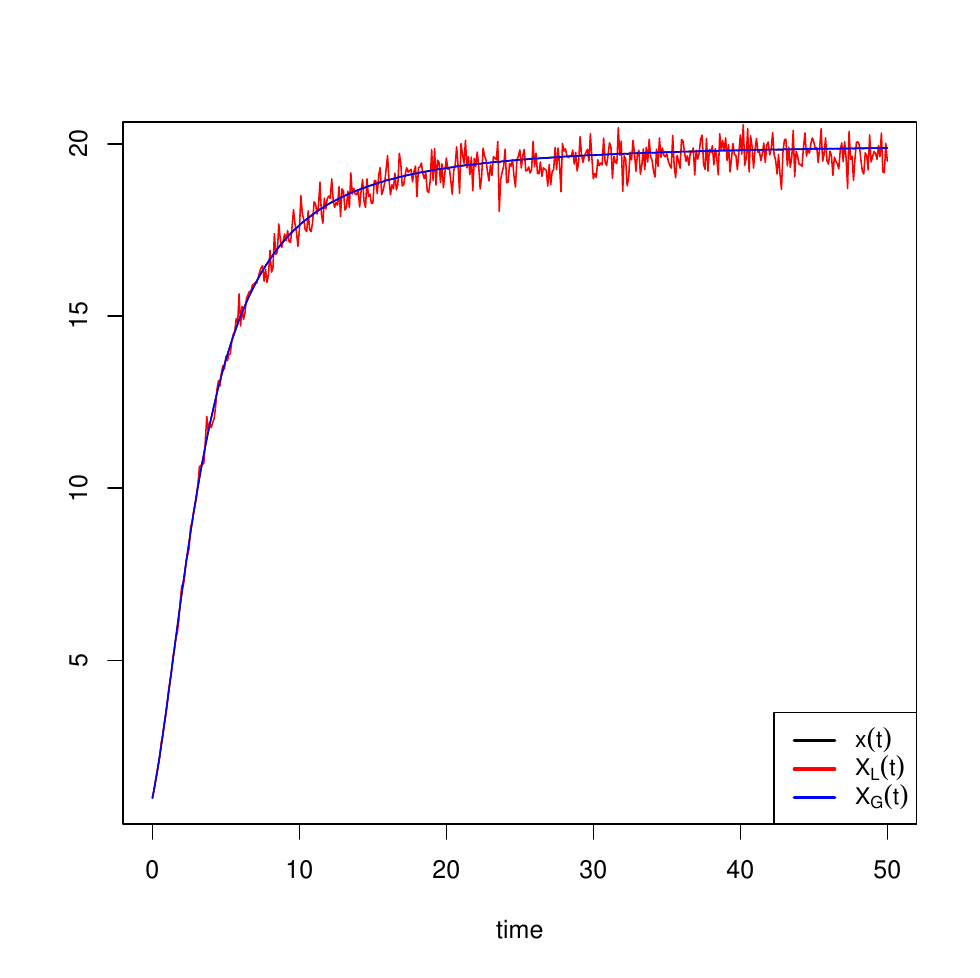}}\quad
\subfigure[$p=1$]
{\includegraphics[width=0.45\textwidth]{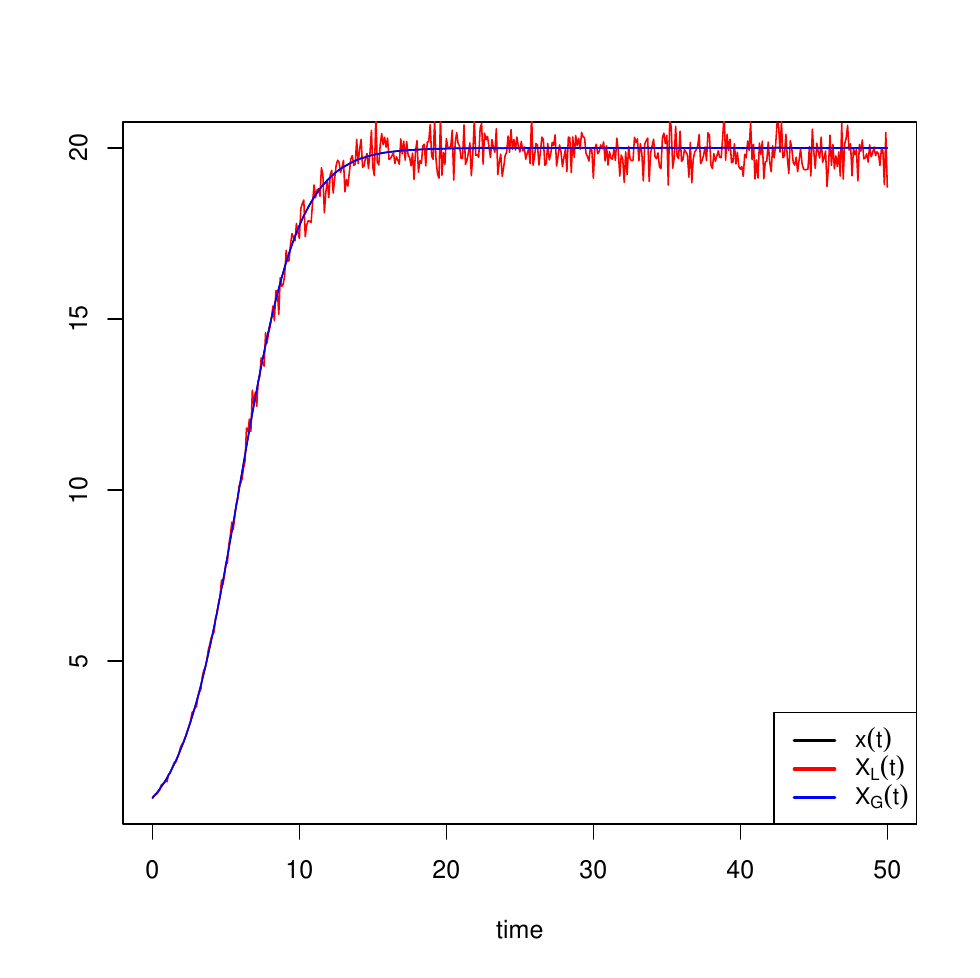}}\\
\subfigure[$p=0.75$]
{\includegraphics[width=0.45\textwidth]{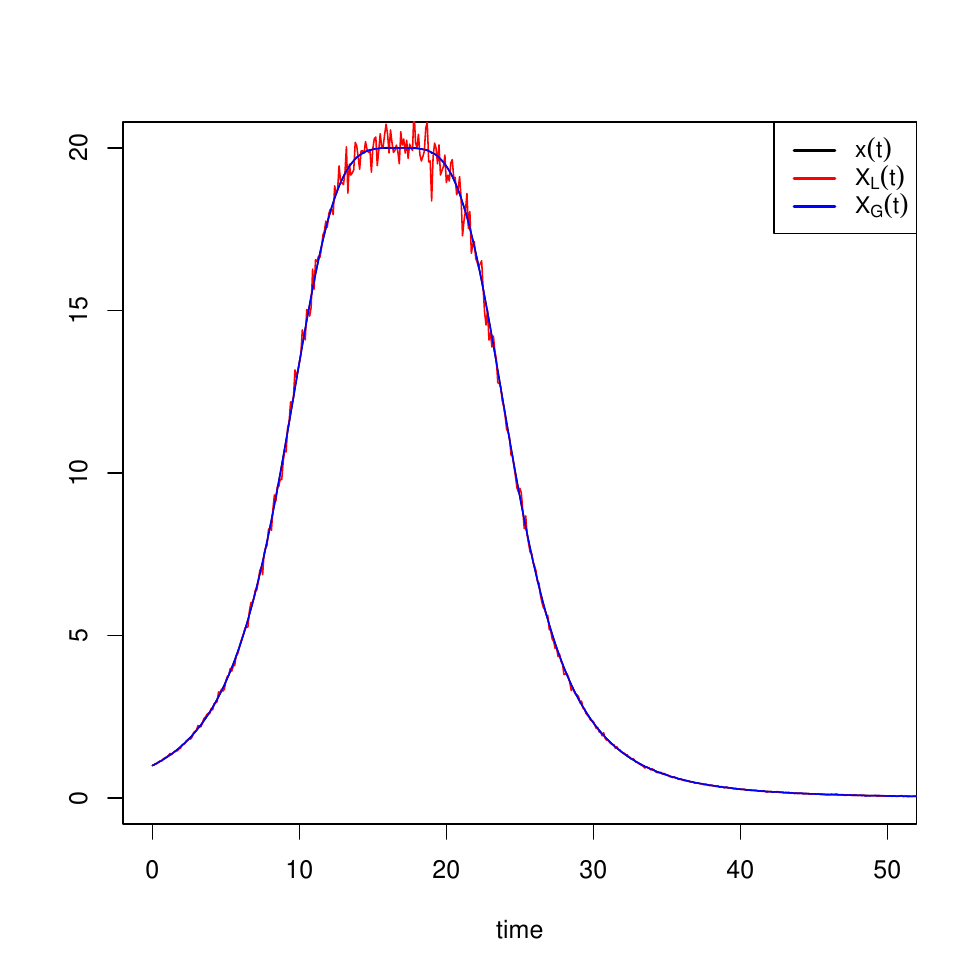}}\quad
\subfigure[$p=2/3$]
{\includegraphics[width=0.45\textwidth]{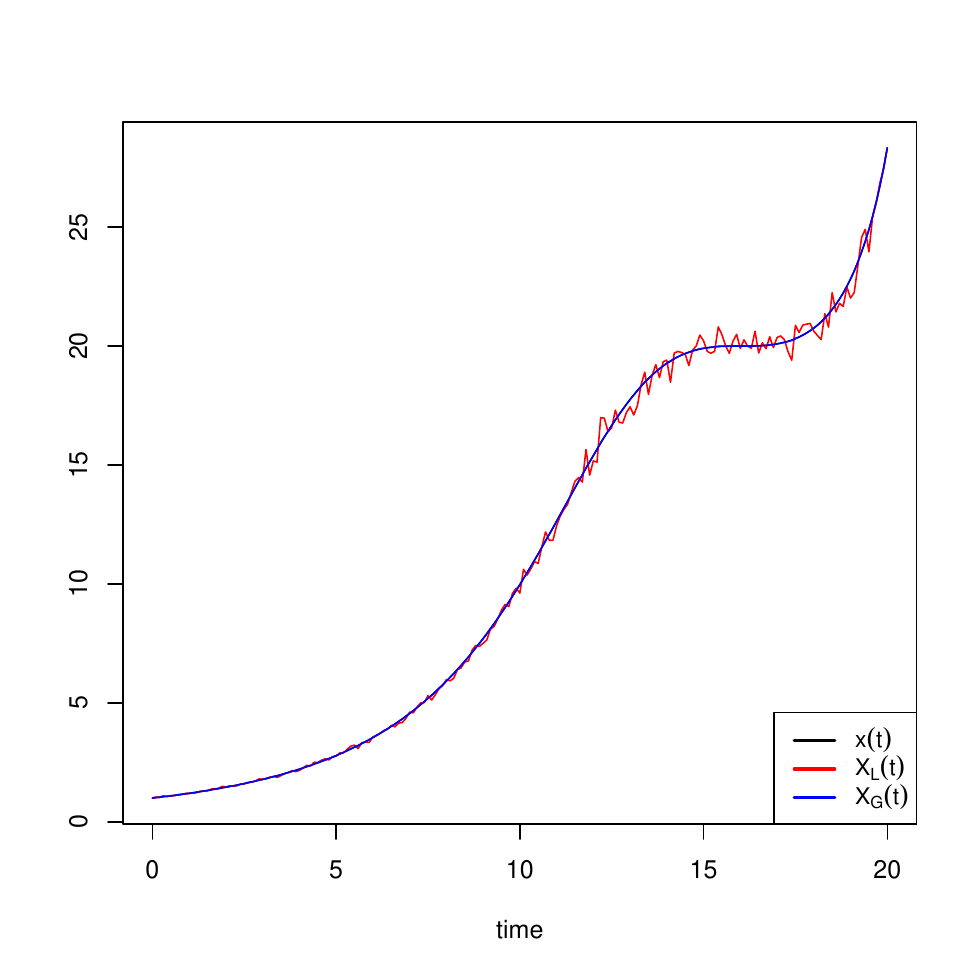}}\\
\subfigure[$p=1/4$]
{\includegraphics[width=0.45\textwidth]{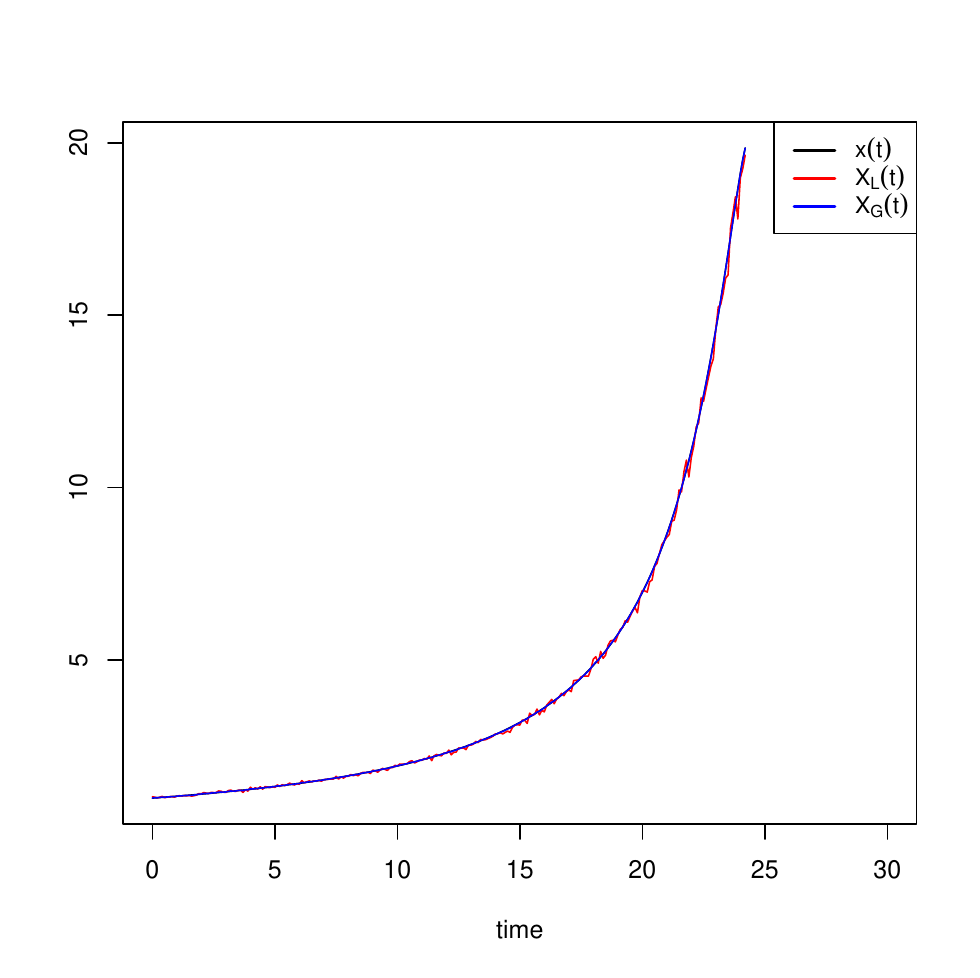}}
\caption{Sample paths of $X_L(t)$ and $X_G(t)$ for $n=1$, $k=20$, $\gamma=0.5$, $x_0=1$, $t_0=0$ and $\sigma=0.02$.}  \label{fig1}
\end{figure}

In this section first we compare the sample paths of the two processes $X_L(t)$ and $X_G(t)$ and then we provide a numerical analysis for the FPT and the FET problems.
As shown in \cite{Albano_2022}, the deterministic curve $x(t)$ in \eqref{ODE} exhibits different behaviors depending on whether $0<p<1$ and $1\leq p<1+\frac{1}{n}$. Clearly, different dynamics are shown also for the processes $X_L(t)$ and $X_G(t)$, since the sample paths of both of them move around to their mean function $x(t)$.
In Figure \ref{fig1} we compare the deterministic curve (black curve) with the sample paths of the processes $X_L(t)$ (red) and $X_G(t)$ (blue) for $n=1, \gamma=0.5, k=20, x_0=1, t_0=0, \sigma=0.02$ and for several choices of $p$ aimed to show the possible dynamics of the processes. Such behaviors are in accordance with those ones discussed in \cite{Albano_2022}. Indeed in Figures \ref{fig1}(a) and \ref{fig1}(b) we have that $1\leq p<1+\frac{1}{n}$ and the paths have a sigmoidal shape going to the carrying capacity $k$; instead, in Figures \ref{fig1}(c), \ref{fig1}(d) and \ref{fig1}(e),  we have chosen $0<p<1$ and the sample paths present a different behavior in each case. In particular, in Figure  \ref{fig1}(c) the ratio $\frac{1}{1-p}$ is even and the trend of sample paths is non monotonic showing an initial increasing followed by a \lq\lq plateau\rq\rq\  around the value $k$, after which the process tends to zero in a decreasing way. In Figure \ref{fig1}(d) $p$ is such that the ratio $\frac{1}{1-p}$ is odd and, after a plateau, the processes indefinitely increase. Finally, in Figure \ref{fig1}(e) the ratio $\frac{1}{1-p}$ is not an integer and the paths reach the value $k$ in a finite time.
In all the figures  the paths of the processes $X_L(t)$ and $X_G(t)$ present oscillations with very different widths. {\color{red} Indeed, the variability of the lognormal process $X_L(t)$ depends on the state in which the process is, while $X_G$ is homoskedastic, therefore it has constant variability. This different variability implies that for the same value of the parameter $\sigma$, the sample paths of $X_G(t)$ have oscillations smaller around the mean function $x(t)$ with respect to the paths of $X_L(t)$. }
For this reason, in the following we choose values of $\sigma$ different for the two processes.
\subsection{Analysis of the process $X_L(t)$}
In our numerical analysis we choose $n=1, \gamma=0.5, k=20, x_0=1, t_0=0$ and we consider several choices of $\sigma$. From  \eqref{boundary} we have that $S(t)=\nu \frac{g(0)}{g(t)}$ with
\begin{equation}\label{gFix}
g(t)=\frac{1}{19}+\Bigl[1-\frac{19^{p-1}}{2} (1-p) t\Bigr]^{1/(1-p)}
\end{equation}
We point out that, although $S(t)$ depends on the parameter $p$, the FPT pdf $g_L[s(t),t| x_0,t_0]$ in \eqref{proportional} does not depend on it. Therefore the different cases analyzed in Figure 1 lead to thresholds having different behaviors but presenting the same FPT pdf.

\begin{figure}[htbp]
\centering
\subfigure[$\sigma=0.02$]{\includegraphics[width=0.45\textwidth]{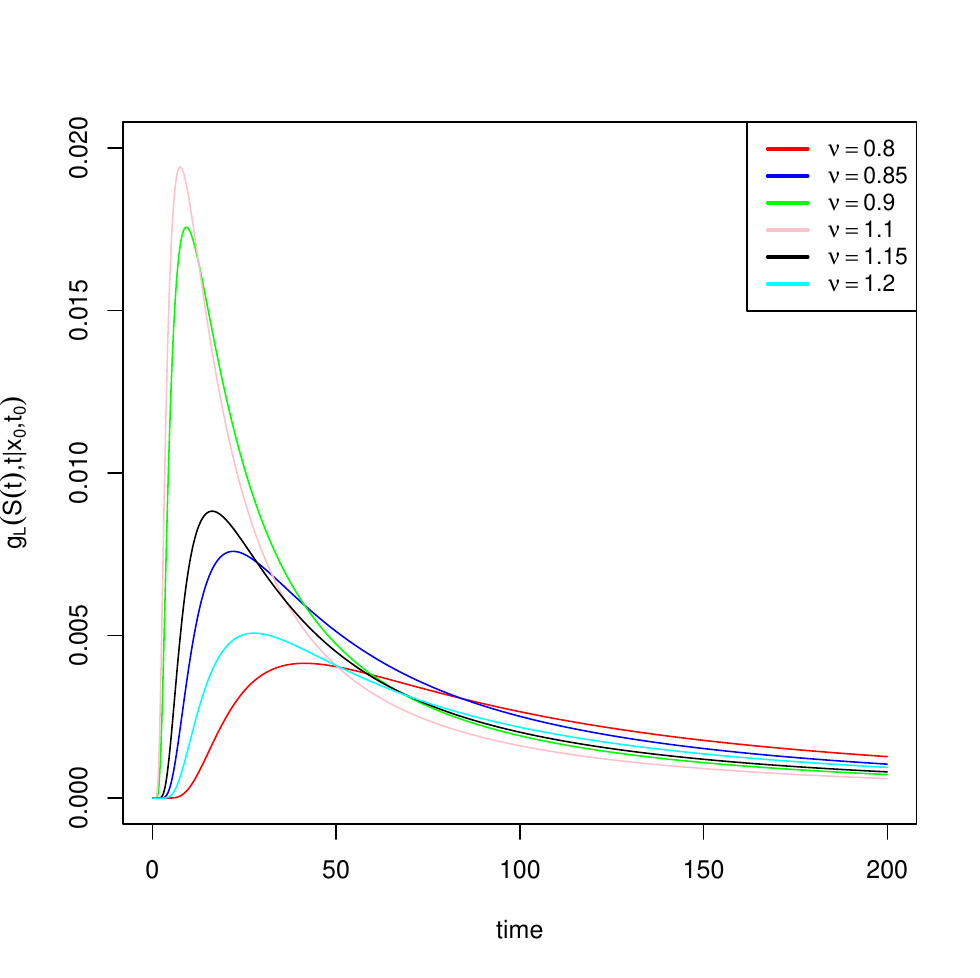}}\quad
\subfigure[$\nu=0.8$]{\includegraphics[width=0.45\textwidth]{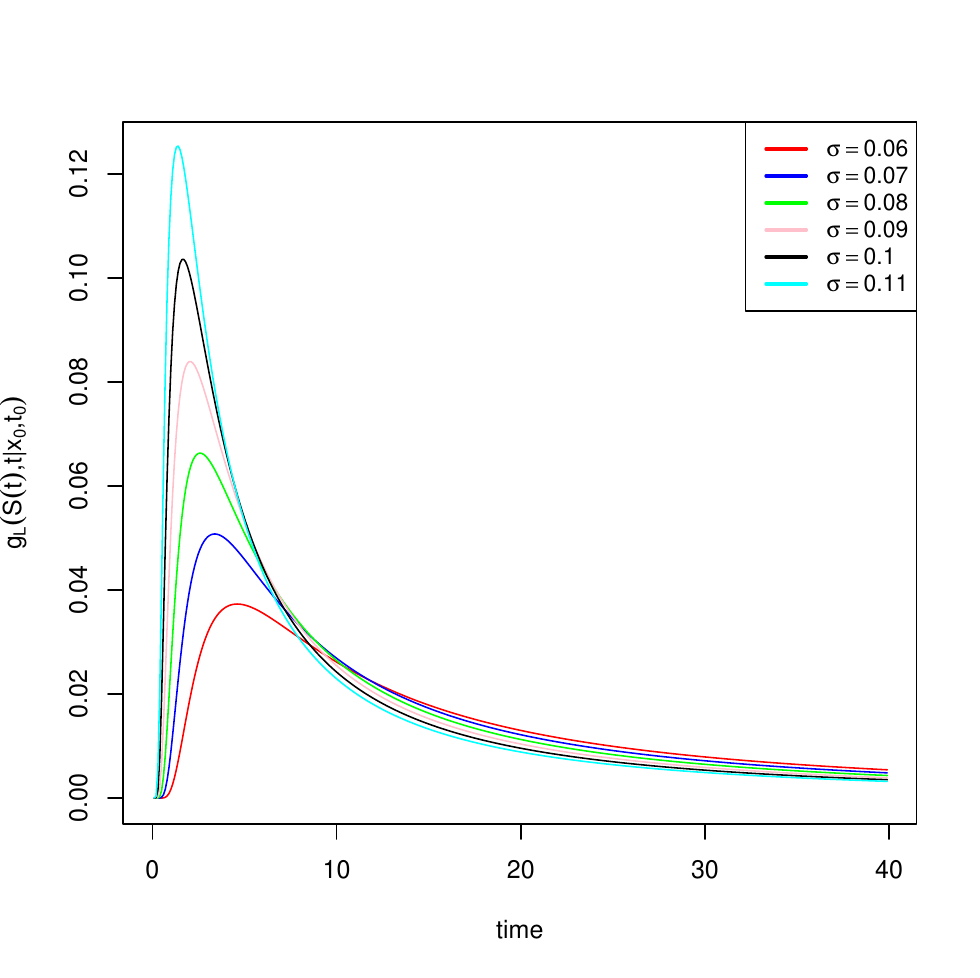}}
\caption{For $n=1, \gamma=0.5, k=20, x_0=1, t_0=0$, the FPT density of $X_L(t)$ for different values of the proportion $\nu$ (on the left) and  for different values of $\sigma$ (on the right).}  \label{fig2}
\end{figure}

In Figure \ref{fig2} we show the FPT pdf of $X_L(t)$ through $s(t)=\nu\frac{g(0)}{g(t)}$ with $g(t)$ given in \eqref{gFix} for different values of the proportion $\nu$ (Figure \ref{fig2}(a)) and  for different values of $\sigma$ (Figure \ref{fig2}(b)). {\color{red} By fixing $\sigma$ (Figure \ref{fig2}(a)), we observe two different behaviors of the FPT pdfs as $\nu$ increases. Precisely, for $0<\nu<1$, the FPT pdf becomes more and more peaked and the maximum is achieved in a shorter time as $\nu$ increases. An inverse behavior is instead shown by the FPT pdf if $\nu\geq 1$. Further, when the proportion $\nu$ is fixed and the width of the oscillations $\sigma$ increases (Figure \ref{fig2}(b)), the abscissa of the maximum of the FPT pdf decreases as $\nu$ increases, while its ordinate increases. This can be verified analytically by studying the instant of time in which the FPT pdf of $X_L(t)$ reaches its maximum, seen as a function of $\sigma^2$ and of $\nu$.

}

Now we consider the FET problem of the process $X_L(t)$ through the region $[s_1(t), s_2(t)]$ with $s_i(t)=\nu_i \frac{g(0)}{g(t)}$ and $g(t)$ defined in \eqref{gFix}.
Also in this case, the FET pdf $\gamma_L(t| x_0,t_0)$ given in \eqref{gamma_percentuali} does not depend on the parameter $p$, while such a dependence is preserved in the boundaries $s_i(t)\ (i=1,2)$.

\begin{figure}[htbp]
\centering
\subfigure[$\nu_2=1.3,\quad \sigma=0.02$]{\includegraphics[width=0.45\textwidth]{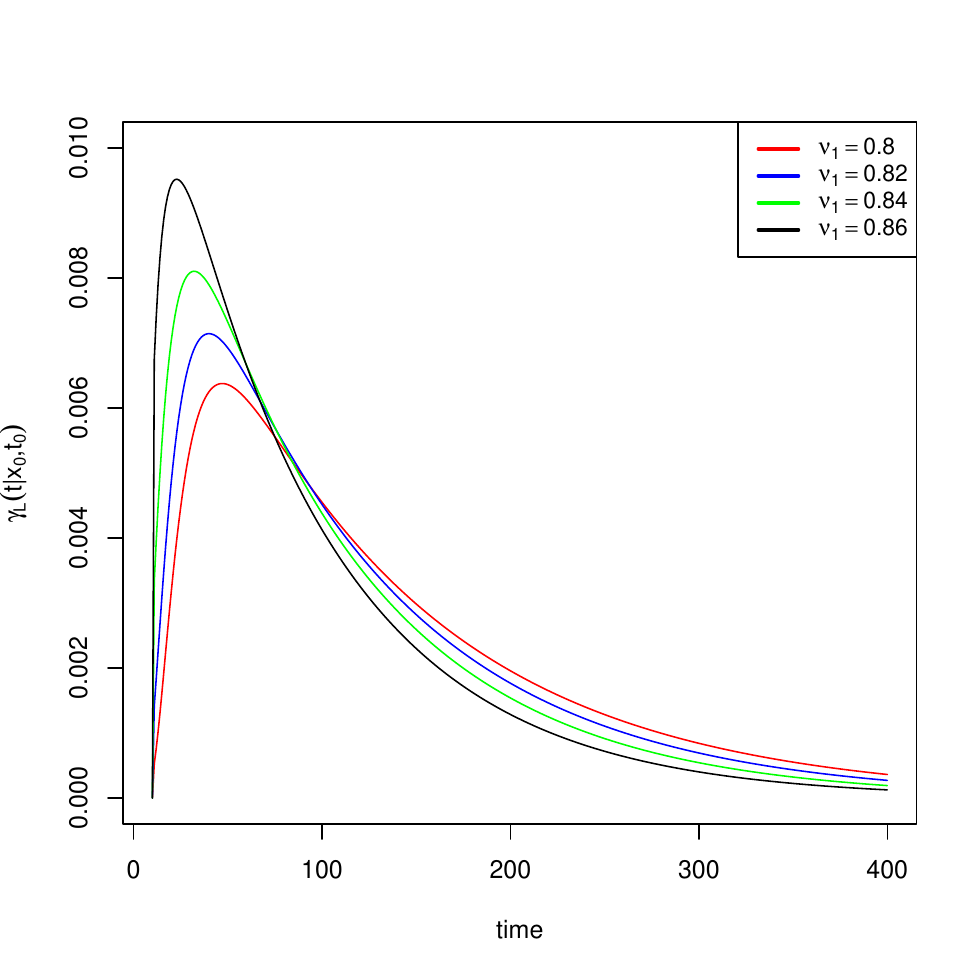}}\quad
\subfigure[$\nu_1=0.8,\quad \sigma=0.02$]{\includegraphics[width=0.45\textwidth]{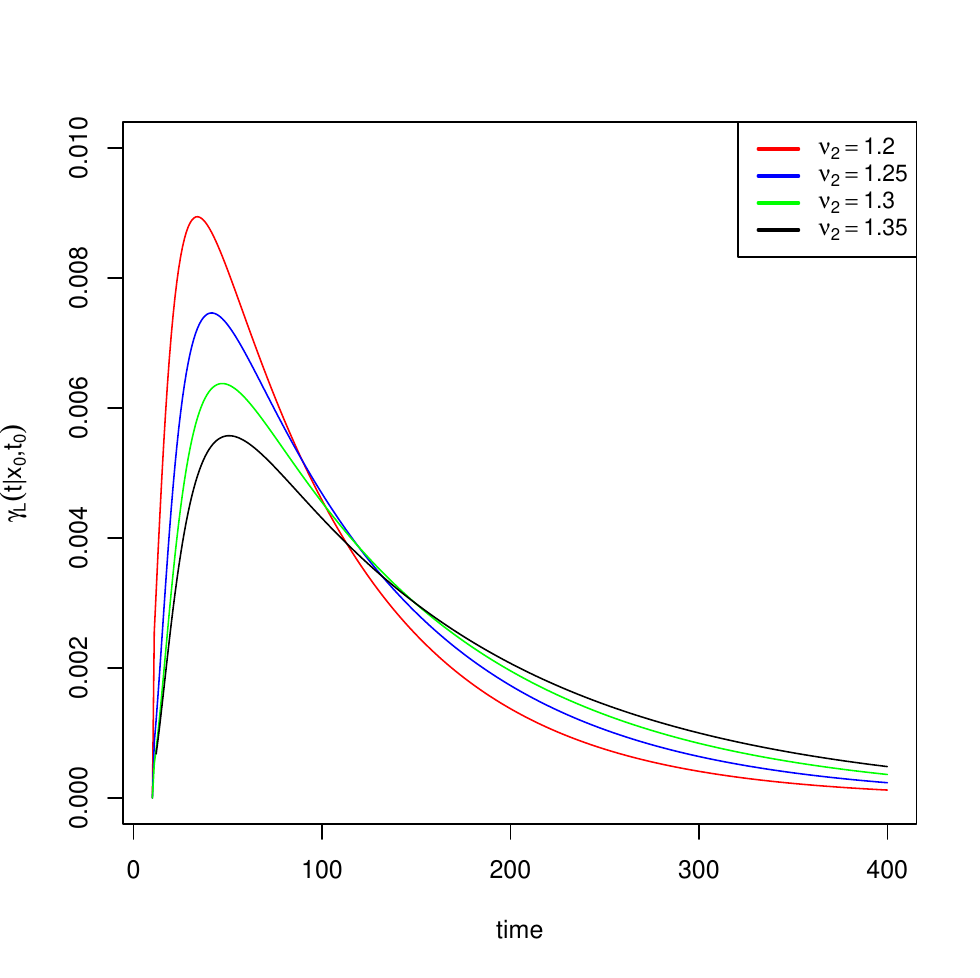}}
\caption{For $n=1, \gamma=0.5, k=20, x_0=1, t_0=0, \sigma=0.02$, the FET pdf of $X_L(t)$ for different values of $\nu_1$ (on the left) and  for different values of $\nu_2$ (on the right).}
\label{fig4}
\end{figure}

\begin{figure}[htbp]
\centering
\subfigure[$\nu_1=0.8,\quad \nu_2=1.2$]{\includegraphics[width=0.45\textwidth]{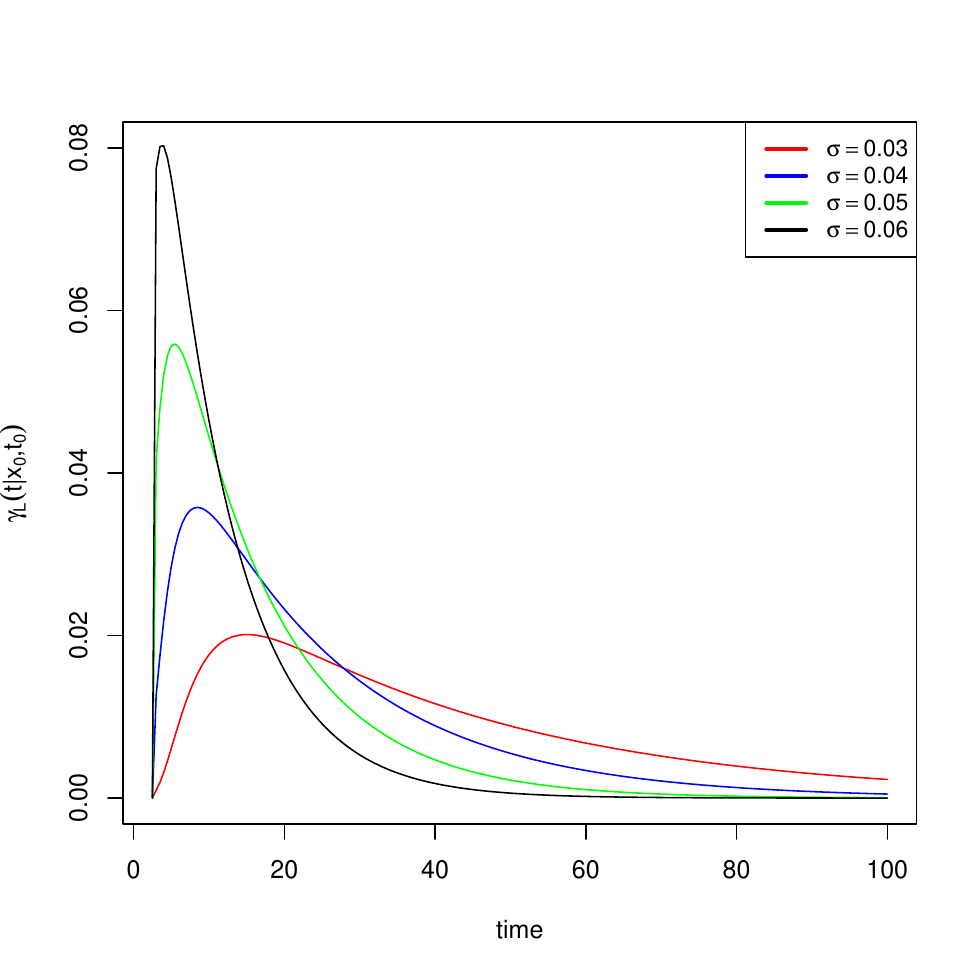}}\quad
\subfigure[$\nu_1=0.8,\quad \nu_2=1.3$]{\includegraphics[width=0.45\textwidth]{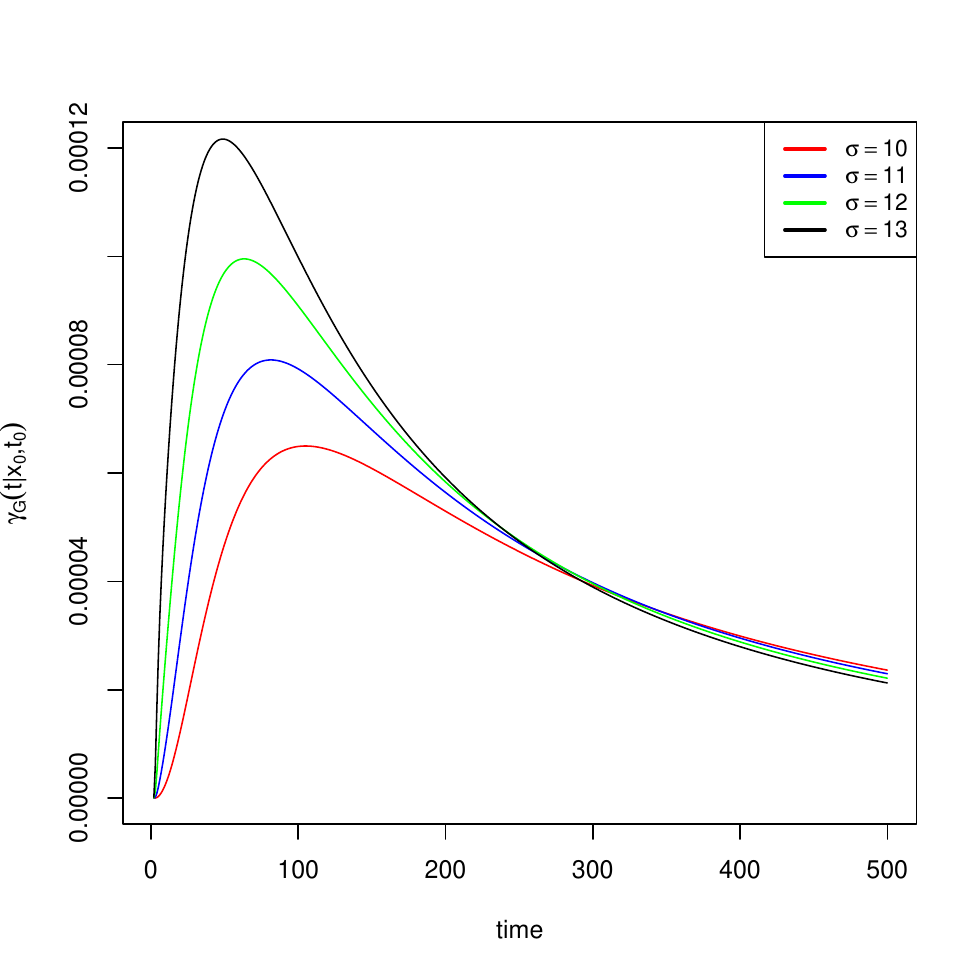}}
\caption{For $n=1, \gamma=0.5, k=20, x_0=1, t_0=0$, the FET pdf of $X_L(t)$ for different values of $\sigma$.}
\label{fig5}
\end{figure}

In Figures \ref{fig4} and \ref{fig5} we have plotted the FET pdf $\gamma_L(t| x_0,t_0)$ for several choices of $\nu_1, \nu_2$ and $\sigma$. In particular in Figure \ref{fig4}(a) we have fixed $\sigma$ and $\nu_2$, we note that the FET pdf becomes more and more peaked as $\nu_1$ increases. This is due to the fact that increasing $\nu_1$ is equivalent to narrowing the region of interest $[s_1(t), s_2(t)]$. Instead, by increasing the parameter $\nu_2$, the region becomes wider and, as shown in Figure \ref{fig4}(b), by fixing $\sigma$ and $\nu_1$, the FET pdf becomes less and less peaked as $\nu_2$ increases.

In Figure \ref{fig5} the FET pdf $\gamma_L(t| x_0,t_0)$ is plotted for several values of $\sigma$, with fixed values of the proportions $\nu_1$ and $\nu_2$. In this case the region is fixed and the amplitude of the oscillations varies. Consequently, the maximum of the FET pdf is higher and it is reached for shorter times. Furthermore, as expected, by enlarging the width of the region $[s_1(t), s_2(t)]$, i.e. comparing Figures \ref{fig5}(a) and \ref{fig5}(b), the maximum of the FET pdf gets lower.

\subsection{Analysis of the process $X_G(t)$}
For the process $X_G(t)$ we consider the FPT pdf $g_G(s(t), t|x_0,t_0)$ given in \ref{gG} and the FET pdf $\gamma_G(t|x_0,t_0)$ in \ref{gammaG}.
As in the case of the lognormal process $X_L(t)$, we choose $n=1, \gamma=0.5, k=20, x_0=1, t_0=0$ and we consider several choices of $\sigma$. From  \eqref{boundary} we have that $S(t)=\nu \frac{g(0)}{g(t)}$ with $g(t)$ given in \ref{gFix}. Also in this case, the FPT pdf $g_G(s(t), t|x_0,t_0)$ and the FET pdf $\gamma_G(t|x_0,t_0)$ do not depend on the parameter $p$ that, instead, affects the boundaries $s(t)$, $s_1(t)$ and $s_2(t)$.

\begin{figure}[htbp]
\centering
\subfigure[$\sigma=0.1$]{\includegraphics[width=0.45\textwidth]{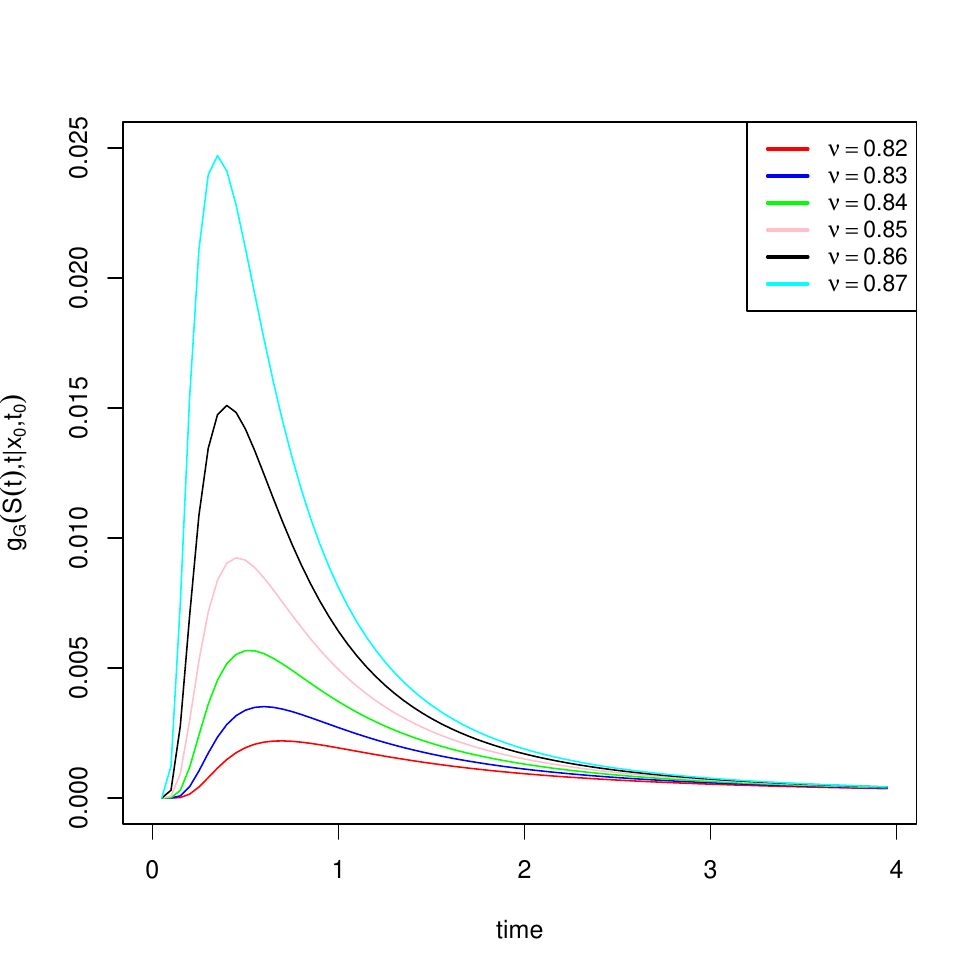}}\quad
\subfigure[$\nu=0.8$]{\includegraphics[width=0.45\textwidth]{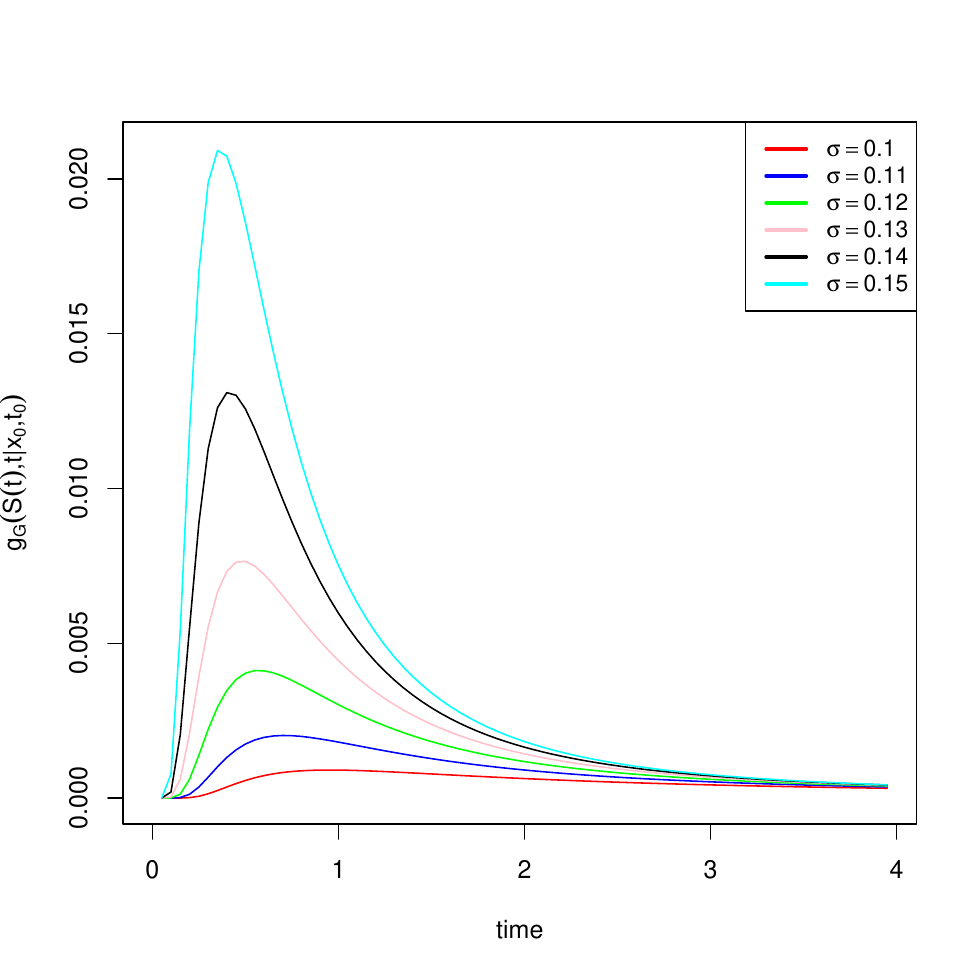}}
\caption{For $n=1, \gamma=0.5, k=20, x_0=1, t_0=0$, the FPT density of $X_G(t)$ for different values of the proportion $\nu$ (on the left) and  for different values of $\sigma$ (on the right).}
\label{fig3}
\end{figure}

In Figure \ref{fig3} the FPT pdf of $X_G(t)$ through $s(t)=\nu\frac{g(0)}{g(t)}$ is plotted for different values of $\nu$ (Figure \ref{fig3}(a)) and  for different values of $\sigma$ (Figure \ref{fig3}(b)). For $\sigma$ fixed (Figure \ref{fig3}(a)), the maximum of the FPT pdf becomes higher and it is reached for shorter times, as $\nu$ increases. Further, the FPT pdf exhibits a similar behavior when the proportion $\nu$ is fixed and the width of the oscillations $\sigma$ increases (Figure \ref{fig3}(b)).

Now we consider the FET problem of the process $X_G(t)$ through the region $[s_1(t), s_2(t)]$ with $s_i(t)=\nu_i \frac{g(0)}{g(t)}$ and $g(t)$ defined in \eqref{gFix}.

\begin{figure}[htbp]
\centering
\subfigure[$\nu_2=1.1,\quad \sigma=10$]{\includegraphics[width=0.45\textwidth]{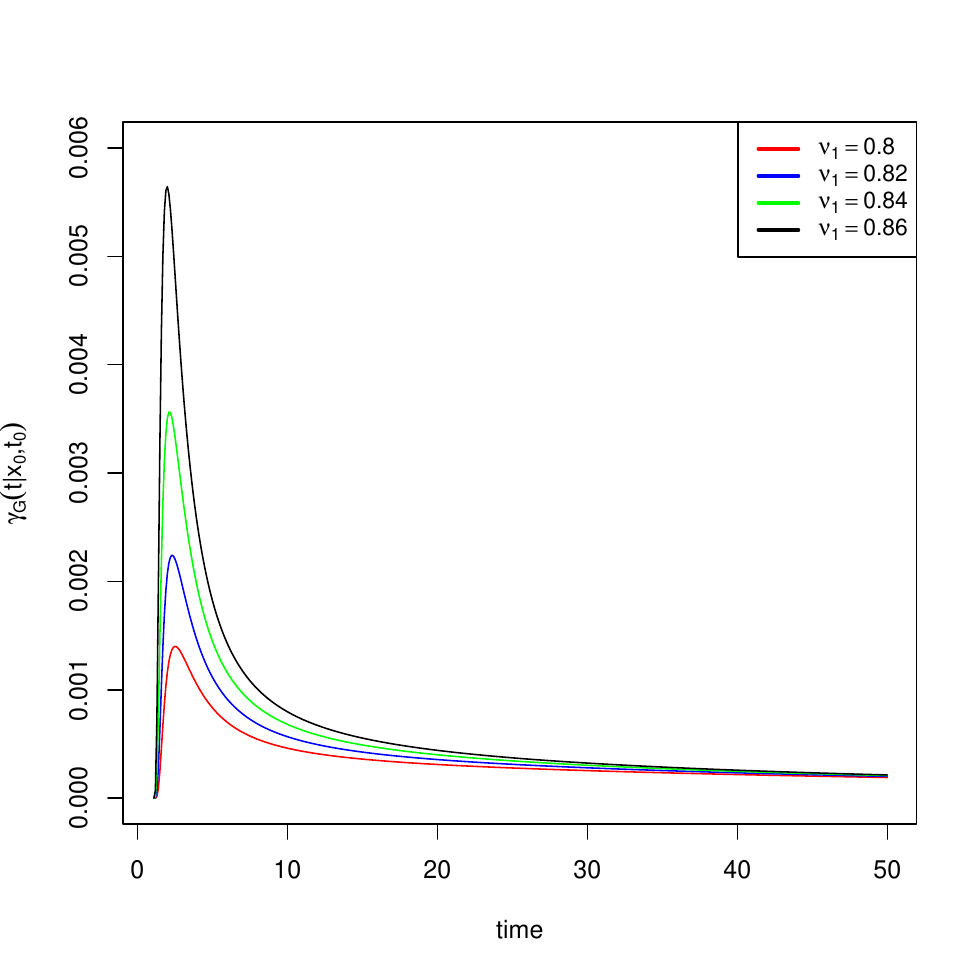}}\quad
\subfigure[$\nu_1=0.8,\quad \sigma=10$]{\includegraphics[width=0.45\textwidth]{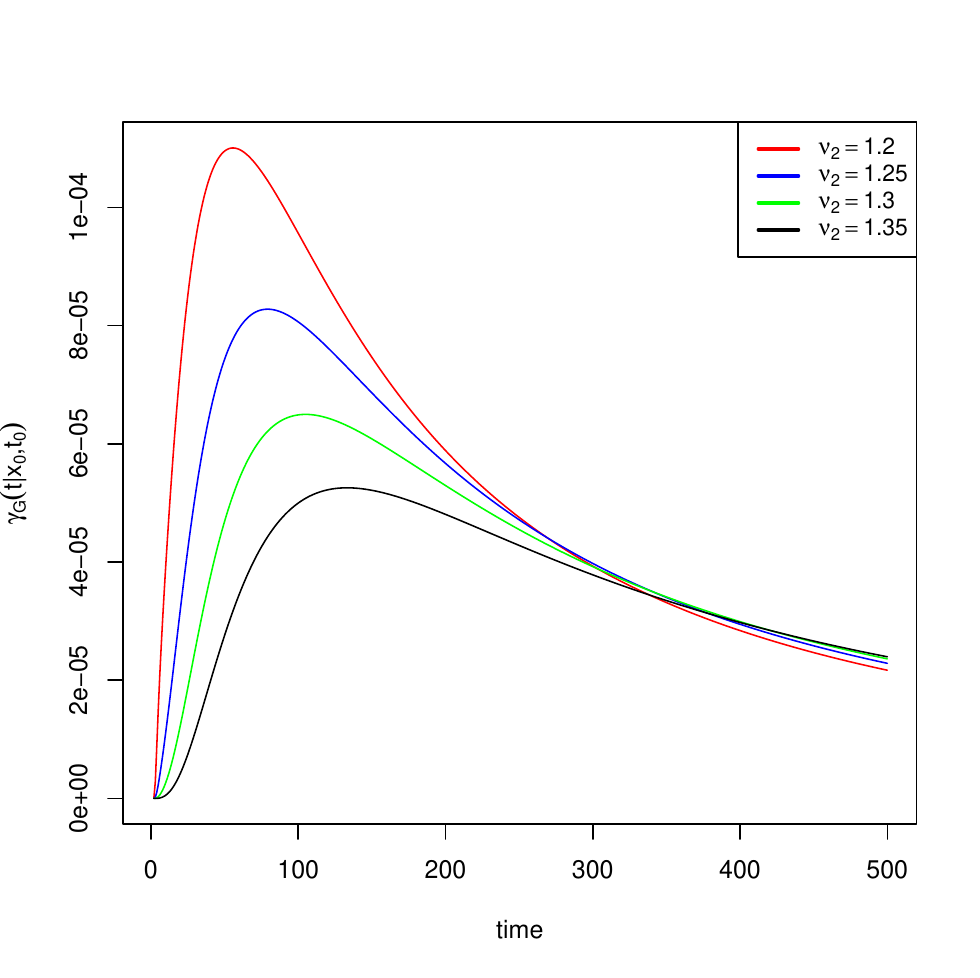}}
\caption{For $n=1, \gamma=0.5, k=20, x_0=1, t_0=0, \sigma=10$, the FET pdf of $X_G(t)$ for different values of $\nu_1$ (on the left) and  for different values of $\nu_2$ (on the right).}  \label{figG}
\end{figure}

\begin{figure}[htbp]
\centering
\subfigure[$\nu_1=0.8,\quad \nu_2=1.2$]{\includegraphics[width=0.45\textwidth]{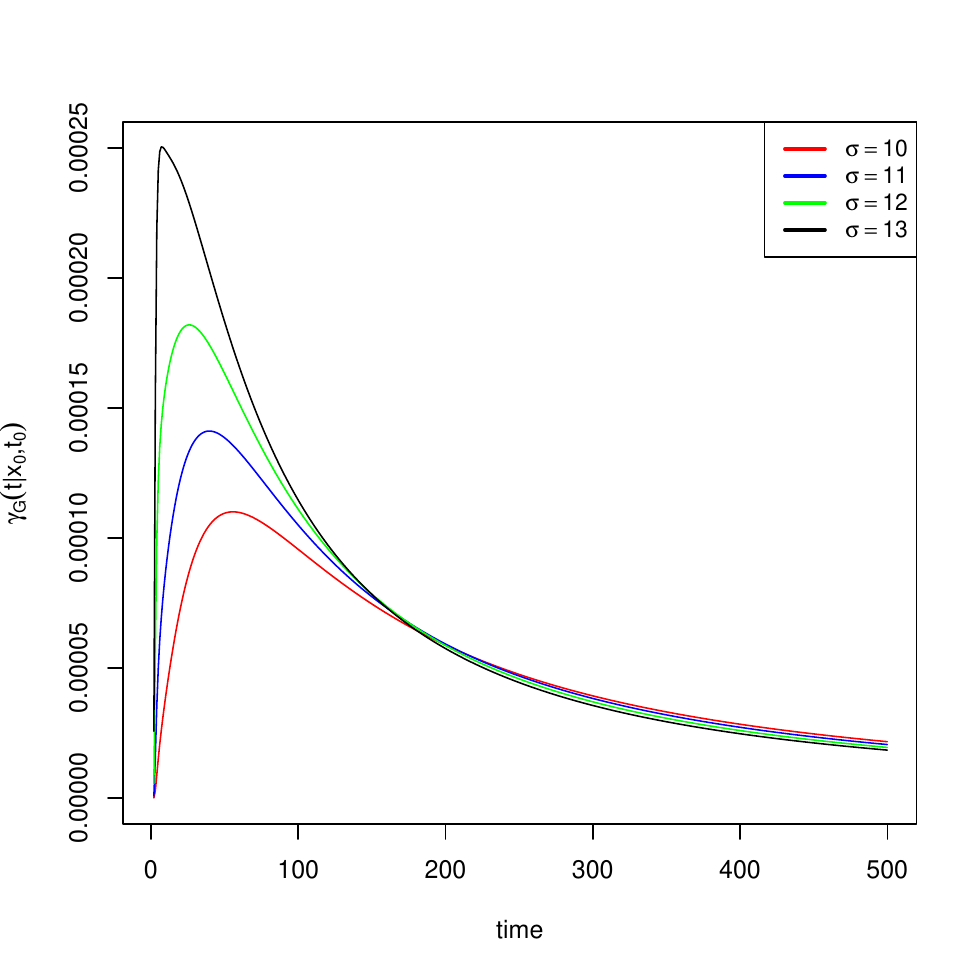}}\quad
\subfigure[$\nu_1=0.8,\quad \nu_2=1.3$]{\includegraphics[width=0.45\textwidth]{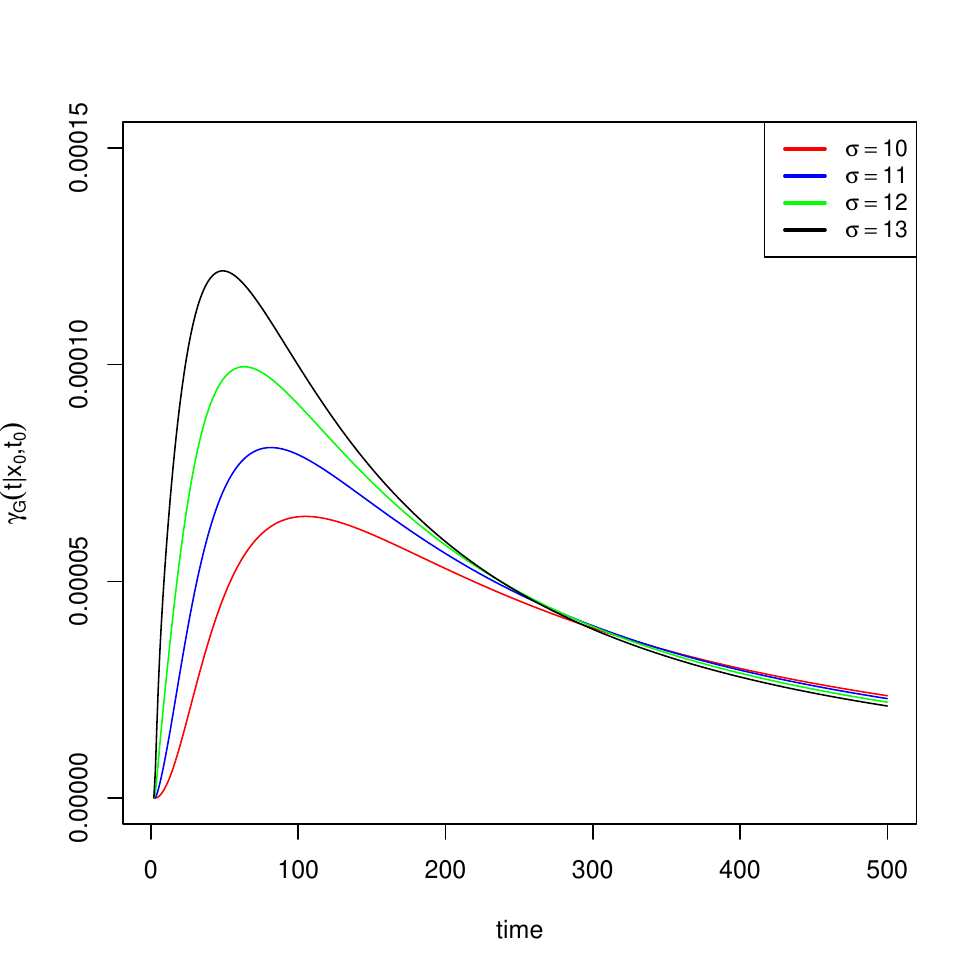}}
\caption{For $n=1, \gamma=0.5, k=20, x_0=1, t_0=0$, the FET pdf of $X_L(t)$ for different values of $\sigma$.}
\label{figGG}
\end{figure}

In Figures \ref{figG} and \ref{figGG} we have plotted the FET pdf $\gamma_G(t| x_0,t_0)$ for several choices of $\nu_1, \nu_2$ and $\sigma$. In particular in Figure \ref{figG}(a) we have fixed $\sigma=10$ and $\nu_2=1.1$, we note that the FET pdf becomes more and more peaked as $\nu_1$ increases, i.e. narrowing the region $[s_1(t), s_2(t)]$. Instead, by increasing the parameter $\nu_2$ and fixing $\sigma$ and $\nu_1$ (Figure \ref{figG}(b)), the FET pdf becomes less and less peaked as $\nu_2$ increases.

In Figure \ref{figGG} the FET pdf $\gamma_G(t| x_0,t_0)$ is shown for several values of $\sigma$, with fixed values of $\nu_1$ and $\nu_2$; as $\sigma$ increases, the maximum of the FET pdf becomes higher and it is reached for shorter times. Finally, comparing the Figures \ref{figGG}(a) and \ref{figGG}(b), we observe that as the width of the region increases, the FET pdf becomes flatter and the tails become heavier.

\section{Concluding remarks}
We have  proposed two stochastic formulations of a general growth equation, including the well known growth equations, such as Malthus, logistic, Bertallanffy, Gompertz.  They have been obtained starting from a new parametrization of the growth equation, by adding an additive and multiplicative noise.
The processes obtained, although having the same mean, which coincides with the deterministic solution of the growth equation, are lognormal Gaussian distributed. This diversity is also evident in the trajectories of the two processes which are differently influenced by the amplitude of the random oscillations.

For these processes we have analyzed  the FPT from a threshold and on the FET from a region delimited by two thresholds. We have provided specific time dependent thresholds, for which there exist closed-forms of the FPT and FET pdf's.
it is interesting to observe that, with appropriate choices of the parameters, the identified thresholds represent percentages of the average of the two processes. Hence, in application contexts, our analysis answers the question how much time is necessary for processes to reach a certain percentage of the average.

Some numerical results, aimed at analyzing the effect of the parameters and thresholds, is finally provided.

\section*{Acknowledgements}
This work was supported in part by the ``Mar\'ia de Maeztu'' Excellence Unit IMAG, reference CEX2020-001105-M,
funded by MCIN/AEI/10.13039/
501100011033/, by the Ministerio de Ciencia e Innovaci\'on, Spain, under Grant PID2020-
1187879GB-100
and by MIUR - PRIN 2017, Project \lq\lq Stochastic Models for Complex Systems\rq\rq. G. Albano and V. Giorno are members of the group GNCS of INdAM.

\end{document}